\documentclass[a4paper]{article}

\usepackage[margin=1in]{geometry}

\usepackage[T1]{fontenc}
\newcommand{\changefont}[3]{
\fontfamily{#1} \fontseries{#2} \fontshape{#3} \selectfont}

\changefont{ptm}{m}{n}

\usepackage{setspace} 
 \usepackage{graphicx}    

\usepackage{color}																																														
\usepackage{epsfig}
\usepackage{amsfonts}
\usepackage{amsmath}
\usepackage{amssymb}
\usepackage{graphics}
\usepackage{mathrsfs}
\usepackage{subfigure}
\usepackage{setspace}
\newcommand \be{\begin{equation}}
\newcommand \ee{\end{equation}}
\newcommand \ba{\begin{eqnarray}}
\newcommand \ea{\end{eqnarray}}

\def\bit{\begin{itemize}}
\def\eit{\end{itemize}}

\usepackage{amsfonts}
\usepackage{amssymb}
\usepackage{graphics}
\usepackage{mathrsfs}
\usepackage{color}
\newtheorem{remark}{Remark}[section]

\newtheorem{theorem}{Theorem}[section]

\newtheorem{example}{Example}[section]

\newtheorem{lemma}{Lemma}[section]

\newtheorem{definition}{Definition}[section]

\long\def\symbolfootnote[#1]#2{\begingroup%
\def\thefootnote{\fnsymbol{footnote}}\footnote[#1]{#2}\endgroup} 

\begin{document}

%

\begin{center}
\Large \textbf{Discontinuous dynamics with grazing points}
\end{center}

\vspace{-0.3cm}
\begin{center}
\normalsize \textbf{M. U. Akhmet$^{a,} \symbolfootnote[1]{Corresponding Author Tel.: +90 312 210 5355,  Fax: +90 312 210 2972, E-mail: marat@metu.edu.tr}$, A. K{\i}v{\i}lc{\i}m$^{a}$} \\
\vspace{0.2cm}
\textit{\textbf{\footnotesize$^a$Department of Mathematics, Middle East Technical University, 06800, Ankara, Turkey}} 
\vspace{0.1cm}
\end{center}

\vspace{0.3cm}

\begin{center}
\textbf{Abstract}
\end{center}

\noindent\ignorespaces

Discontinuous dynamical systems with grazing solutions are discussed. The group property, continuation of solutions, continuity and smoothness of motions are thoroughly analyzed. A variational system around a grazing solution which depends on near solutions is constructed. Orbital stability of  grazing cycles is  examined by linearization. Small parameter method is extended for analysis of neighborhoods of grazing orbits, and grazing bifurcation of cycles is observed in an example. Linearization around an equilibrium grazing point is discussed. The mathematical background of the study relies on the theory of discontinuous dynamical systems \cite{ref1}. Our approach is analogous to that one of the continuous dynamics analysis and results can be extended on functional differential, partial differential equations and others. Appropriate illustrations with grazing limit cycles and bifurcations are depicted to support the theoretical results.

\vspace{0.2cm}
 
\noindent\ignorespaces \textbf{Keywords:} Discontinuous dynamical systems; Grazing points and orbits; Axial and non-axial grazing; Variational system; Orbital stability; Small parameter; Bifurcation of cycles; Impact mechanisms

\section{Introduction} \label{secintro}

Vibro-impacting systems provide  examples of non-linear dynamical systems, exhibiting new levels of complicated dynamics due to their non-smoothness.  Grazing phenomenon is one of the  attractive  features  for  these  dynamics \cite{nusse94}-\cite{Nord1}, \cite{Luo2005a}-\cite{luotime}. There are two approaches for the definition of grazing in literature.  One is presented in the studies of   Bernardo,   Budd and   Champneys  \cite{budd2001,budd1998}, Bernardo and Hogan \cite{Bernardo-Hogan2010}, and Luo \cite{Luo2005a}-\cite{luotime}.  In these studies,  it is asserted that grazing occurs when a trajectory hits the surface of discontinuity tangentially.    In  \cite{Nord1}-\cite{Nord2006}, Nordmark defines grazing   as the approach of the velocity to zero in the neighborhood of the surface of discontinuity which is the case of the studies conducted in  \cite{budd2001}-\cite{Bernardo-Hogan2010}, \cite{Luo2005a}-\cite{luotime}. Our comprehension of grazing in this paper is   close to  that  one in   \cite{budd1998},\cite{Bernardo-Hogan2010},\cite{Luo2005a}.  In the paper \cite{donde-hiskens}, grazing   is considered as a bounding case which separates regions of quite different dynamic behaviors. It is understood that the system trajectory makes tangential contact with an event triggering hypersurface. The shooting, continuation and optimization methods are developed and illustrated for both transient and grazing phenomena. It is exemplified by utilizing power electronics and robotics. In \cite{Piir-virgin-Champhneys}, the grazing periodic orbit and its linearization are obtained by means of a numerical continuation method for hybrid systems. Applying this, the normal-form coefficients are evaluated, which in this case imply the occurrences a jump to chaos and period-adding cascade. The necessary and sufficient conditions of the general discontinuous boundary are expounded in \cite{Luo2006}. In \cite{luo22}, by means of non-stick mapping, the necessary and sufficient conditions for the grazing in periodically forced linear oscillator with dry friction are obtained.  By constructing special maps such as zero time discontinuity mapping \cite{ref24}-\cite{ref38} and Nordmark map \cite{Nord1}-\cite{Nord2006}, the existence of periodic solution and their stability  were   investigated  in  mechanical systems.

In this  paper, we model the dynamics with grazing impacts by utilizing differential equations with impulses at variable moments and  applying  the   methods   of \cite{ref1}-\cite{MAkhmethod}.  As a consequence of such methods,  the role of the mappings \cite{Nord1}-\cite{Nord2006}  is diminished. One   can observe  that   a trajectory   at   a grazing point    may    have tangency    to   the   surface of discontinuity,  which  is  parallel to   one   or several coordinate    axises. Particularly  it  means the   velocity   approaches to  zero  \cite{Nord1}-\cite{Nord2006}. Then, we   will say  about the axial   grazing.  Otherwise, grazing  is non-axial.   This research contains the analysis of both axial and non-axial grazing.

In \cite{ref47}, it is observed through simulations and experiments that the coefficient of restitution  depends on the impact velocity of the particle by considering both the viscoelastic and the plastic deformations of particles occurring at low and high velocities, respectively. It has been proposed in \cite{ref31} that at low impact velocities and for most materials with linear elastic range, the coefficient of restitution   is of the   form $R(v) = 1-av,$ where $v$  is the velocity before collision and $a$ is a  constant.  Also, for low impact velocities the restitution law can be considered  quadratic \cite{repeated}. This   is why,  we  will   use  non-constant   restitution  coefficients in   models   with   impacts   in this   study.

For     investigation  of  autonomous differential equations, it is convenient  to utilize properties of  dynamical systems. They are the  group property, continuation of solutions in both time directions, continuity and differentiability in  parameters.  The studies of discontinuous dynamical systems  with transversal  intersections of   orbits  and surfaces, $B-$ smooth discontinuous flows,  can be   found  in \cite{ref1,marat2005}.   In this research,  the dynamics   is   approved for systems with  grazing orbits.  Moreover, the definitions of orbital stability and asymptotic phase are adapted to grazing cycles. The orbital stability theorem is proved, which can not  be underestimated   for   theory   of impact  mechanisms.


The remaining part of the paper is organized as follows. In the next section, we will introduce necessary notations, definitions and  theorems to specify discontinuous dynamical systems. In   Section 3,  it    is shown   how  the   dynamics  can be   linearized   around grazing   orbits.   In Section 4, the theorem of orbital stability is adapted for grazing cycles.  In Section 5, the small parameter analysis is   applied  near grazing orbits  and bifurcation   of  cycles  is observed. In last  three  sections, examples are presented  to actualize the theoretical results numerically and analytically. Finally, Conclusion  covers a summary of our study. 
\section{Discontinuous dynamical systems} \label{secdss}
%

 Let  $\mathbb{R},$ $\mathbb{N}$ and $\mathbb{Z}$ be the sets of all real numbers, natural numbers and integers, respectively. Consider the set $D\in\mathbb{R}^n$ such that $D=\cup D_i,$ where $D_i,$ $i=1,2,\ldots,k,$  components  of $D,$ are disjoint open connected subsets of $\mathbb{R}^n.$ To describe the surface of discontinuity, we present a two times continuously differentiable function $\Phi: D^{r}\rightarrow \mathbb{R}^n.$ The set can be defined as $\Gamma=\Phi^{-1}(0)$ and is a closed subset of $\bar {D},$  where $\bar {D}$ is the closure of $D.$  Denote $\partial\Gamma$ as the boundary of $\Gamma.$  One can easily  see   that  $\Gamma=\cup_{i=1}^{k} \Gamma_i$, where $\Gamma_i$ are parts  of the surface of discontinuity in the components of $D.$ Denote $\tilde{\Gamma}=J(\Gamma),$  $\tilde{\Phi}(x)=\Phi (J^{-1}(x)).$  Denote   an $r-$ neighborhood of $D$ in $\mathbb{R}^n$ for a fixed $r>0$ as $D^{r}.$ Let $\Gamma^{r}$ be the $r-$ neighborhood of $\Gamma$ in $\mathbb{R}^n,$ for a fixed $r>0$  and define  functions $J:\Gamma^r\rightarrow D^{r}$ and $\tilde{J}:\tilde{\Gamma}^r\rightarrow D^{r},$ such that, $J(\Gamma), \tilde{J} (\tilde{\Gamma})\subset D.$ Assume that a function $f(x):D^{r}\rightarrow \mathbb{R}^n$ is continuously differentiable in $D^{r}.$ Set the gradient vector of $\Phi$ as $\nabla\Phi(x).$ 

The following definitions will be utilized in the remaining part of the paper. Let $x(t-)$ be the left limit position of the trajectory and $x(t+)$ be the right limit of the position of the trajectory at the moment $t.$ Define  $\displaystyle{\Delta x(t):= x(t+)-x(t-)}$  as the jump operator for a function $x(t)$ such that $x(t)\in\Gamma$ and $t$ is a moment of discontinuity (discontinuity moment). In other words, the discontinuity moment  $t$ is the moment when the trajectory meets the surface of discontinuity $\Gamma.$   The function $I(x)$ will be used in the following part of the paper which is defined as  $I(x):=J(x)-x,$ for $x\in\Gamma.$

The following assumptions  are  needed  throughout this paper.
\begin{itemize}
\item[(C1)] $\nabla\Phi(x)\neq 0$ for all $x\in \Gamma,$ 
\item[(C2)] $J\in C^{1}(\Gamma^r )$ and $\det\Big[\frac{\partial J(x)}{\partial x}\Big]\neq 0,$ for all $x\in \Gamma^r\setminus \partial\Gamma,$
\item[(C3)] $\Gamma\bigcap \tilde{\Gamma}\subseteq \partial\Gamma\cap \tilde{\partial\Gamma},$ 
\item[(C4)]  $\langle \nabla\Phi(x),f(x)\rangle \neq0 $ if $x\in\Gamma\setminus\partial\Gamma,$
\item[(C5)]  $\langle \nabla\tilde{\Phi}(x),f(x)\rangle \neq0 $ if $x\in\tilde{\Gamma}\setminus\partial\tilde{\Gamma},$
\item[(C6)] $J(x)=x$ for all $x\in\partial\Gamma,$ 
\item[(C7)]  $\tilde{J}(x)=x$ for all $x\in\partial\tilde{\Gamma}.$ 
 
\end{itemize}
  One can verify that $\tilde{\Gamma}=\{x\in D| \tilde{\Phi}(x)=0\}$ and $\tilde{J}(x)\neq x$ on $\tilde{\Gamma}$ since of $(C2).$ Condition $(C1)$ implies that for every $x_0 \in \Gamma,$ there exist a number $j$ and a function $\phi_{x_0}(x_1, \ldots,x_{j-1},x_{j+1},\ldots,x_n)$ such that $\Gamma$ is the graph of the function $x_j =\phi_{x_0}(x_1, \ldots ,x_{j-1},x_{j+1},\ldots,x_n)$ in a neighborhood of $x_0.$ Same is true for every $x_0\in\tilde{\Gamma}.$ Moreover, $\nabla\tilde{\Phi}(x)\neq 0,$ for all $x\in\tilde{\Gamma},$ can be verified by using the condition $(C2).$ The conditions $(C2),(C6),(C7),$ imply that the equality $\tilde{J}(x)=x,$ is true  for all $x\in\partial\Tilde{\Gamma}.$

   Let $\mathscr{A}$ be an interval in $\mathbb{Z}.$ We say that  the   strictly  ordered set $\theta=\{\theta_i\}, i \in \mathscr{A},$ is a \textit{$B-$sequence } \cite{ref1} if one of the following alternatives holds: $(i)$ $\theta=\emptyset,$ $(ii)$ $\theta$ is a nonempty and finite set, $(iii)$ $\theta$ is an infinite set such that $|\theta_i|\rightarrow\infty$ as $i\rightarrow\infty.$ In what follows, $\theta$ is assumed to be a $B-$sequence  .

The main object of our discussion is the following system,
\begin{equation}\label{eq:graz1}
\begin{aligned}
&x'=f(x), \\
&\displaystyle{\Delta x|_{x\in \Gamma}}= I(x).
\end{aligned}
\end{equation}

In order to define a solution of (\ref{eq:graz1}), we need the following  function and spaces.

A function $\phi(t):\mathbb{R}\rightarrow\mathbb{R}^{n}, \ n\in\mathbb{N}, \,\theta $  is a  $B-$sequence, is from the set $PC(\mathbb{R},\theta)$ if it : $(i)$  is left continuous, $(ii)$ is continuous, except, possibly, points of $\theta,$ where it has discontinuities of the first kind. 

A function $\phi(t)$ is from the set $PC^{1}(\mathbb{R},\theta)$ if $\phi(t),\phi^{\prime}(t)\in PC(\mathbb{R},\theta),$ where the derivative at points of $\theta$ is assumed to be the left derivative. If $\phi(t)$ is a solution of (\ref{eq:graz1}), then it is required that it belongs to $PC^1(\mathbb{R},\theta)$ \cite{ref1}.

We say   that   $x(t):\mathscr{I}\rightarrow \mathbb{R}^n, \mathscr{I} \subset \mathbb R,$ is a solution of (\ref{eq:graz1})   on  $\mathscr{I}$ if    there   exists   an extension  $\tilde x(t)$ of  the function   on  $\mathbb R$   such that   $\tilde x(t)\in PC^1(\mathbb{R},\theta),$   the equality $x'(t)=f(x(t)), \ t\in \mathscr{I},$   is true  if $x(t)\notin\Gamma,$  $x(\theta_i+)=J(x(\theta_i))$   for  $x(\theta_i)\in\Gamma$ and $ x(\theta_i+)\in\tilde{\Gamma},$  $\theta_i \in \mathscr{I}.$  If  $\theta_i$ is a discontinuity moment of $x(t),$   then  $x(\theta_i)\in \Gamma,$ for $\theta_i>0$ and $x(\theta_i)\in \tilde{\Gamma},$ for $\theta_i<0.$ If  $x(\theta_i)\in\partial\Gamma$ or $x(\theta_i)\in\partial\tilde{\Gamma},$ then $x(\theta_i)$ is a point of discontinuity   with zero jump. 

\begin{definition}\label{grazingpoint} A point $x^*$ from $\partial{\Gamma}$ or $\partial{\tilde{\Gamma}}$ is a \textit{grazing point} of system \eqref{eq:graz1} if  $\langle \nabla\Phi(x^*),f(x^*)\rangle=0 $ or\\ $\langle \nabla\tilde{\Phi}(x^*),f(x^*)\rangle=0,$ respectively. If at least one of coordinates of  $\nabla\tilde{\Phi}(x^*)$  is zero  then the grazing is \textit{axial}, otherwise it is \textit{non-axial}.
\end{definition}
\begin{definition}\label{grazingorbit} An orbit $\gamma(x^*)=\{x(t,0,x^*)| x^*\in D, \ t\in \mathbb{R}\}$ of (\ref{eq:graz1}) is grazing if there exists at least one grazing point on the orbit.
\end{definition}



Consider a solution $x(t):\mathbb{R}\rightarrow\mathbb{R}^{n}$ and $\{\theta_i\}$ be the moments of the discontinuity, they are the moments where solution $x(t)$  intersects $\Gamma$ as time increases and  the moments when the solution it intersects $\tilde{\Gamma}$  as time decreases.

A solution $x(t)=x(t,0,x_0),$ $x_0 \in D$ of (\ref{eq:graz1}) locally exists and is unique if the conditions $(C1)-(C3)$ are valid \cite{ref1}.

In what follows, let $\|\cdot\|$ be the Euclidean norm, that is for a vector $x=(x_1,x_2,\ldots,x_n)$ in $\mathbb{R}^n,$ the norm is equal to $\sqrt{x_1^2+x_2^2+\ldots +x_n^2}.$

The following condition for (\ref{eq:graz1}) guarantees that any set of discontinuity moments  of the system  constitutes a $B-$ sequence and we call the condition  \textit{$B-$ sequence condition}. 
\begin{itemize}
\item[(C8)] $\sup_D \|f(x)\|<+\infty,$ and  $\inf_{x_0\in\tilde{\Gamma}}(x_0,y(\zeta,0,x_0) )>0.$ 
\end{itemize}
In \cite{ref1}, some other   $B-$ sequence  conditions are   provided.

We will request for discontinuous dynamical systems that any sequence of discontinuity moments to be a $B-$ sequence.

Let us set the system 
\begin{equation}\label{ode}
y'=f(y)
\end{equation}
for the possible usage in the remaining part of the paper.

Consider a solution $y(t,0,x_0),$ $x_0\in\tilde{\Gamma},$ of  \eqref{ode}. Denote the first meeting point  of the   solution  with the surface $\Gamma,$ provided   the point  exists, by $y(\zeta,0,x_0).$
The following conditions are sufficient for the continuation property. 
\begin{itemize}
\item [(C9)](a) Every solution $y(t,0,x_0),$ $x_0\in D,$ of \eqref{ode} is continuable to either $\infty$ or $\Gamma$ as time increases, \\
(b) Every solution $y(t,0,x_0),$ $x_0\in D,$ of \eqref{ode} is continuable to either $-\infty$ or $\tilde{\Gamma}$ as time decreases.\end{itemize}

To verify the continuation of the solutions of (\ref{eq:graz1}), the following theorems can be applied.

\begin{theorem}\label{continue1}\cite{ref1}   Assume that  conditions $(C8)$ and $(C9)$  are valid. Then, every solution $x(t)=x(t,0,x_0 ),$ $x_0 \in D$ of (\ref{eq:graz1}) is continuable on $\mathbb{R}.$
\end{theorem}

Now, we will present a condition which  is sufficient for the \textit{group property}. 

\begin{itemize}
\item[(C10)] For all $x_0\in D,$ the solution $y(t,0,x_0)$ of \eqref{ode} does not intersect $\tilde{\Gamma} $ before it meets the surface $\Gamma$ as time increases.
\end{itemize}
In other words, for each $x_0\in D$ and a positive number $s$ such that $y(s,0,x_0)\in\tilde{\Gamma},$ there exists a number $r,$ $0\leq r<s,$ such that $y(r,0,x_0)\in\Gamma.$

 It is easy to verify that the condition $(C10)$  is equivalent to the assertion that for all $x_0\in D,$ the solution $y(t,0,x_0)$ of \eqref{ode} does not intersect $ \Gamma$ before it meets the surface $\tilde{\Gamma}$ as time decreases. In other words, for each $x_0\in D$ and a negative number $s$ such that $y(s,0,x_0)\in\tilde{\Gamma},$ there exists a number $r,$ $s<r\leq 0,$ such that $y(r,0,x_0)\in\tilde{\Gamma}.$



\begin{theorem}(The group property)
Assume that conditions (C1)-(C10) are  valid. Then, \\ $x(t_2,0,x(t_1, 0,x_0))=x(t_2 +t_1 ,0,x_0),$ for all $t_1 , t_2 \in\mathbb{R}.$
\end{theorem}

\noindent\textbf{Proof.}  Denote by $\xi(t)=x(t+\bar{t})$, for a fixed $\bar{t}\in\mathbb{R}.$  It can be verified that the sequence $\{\theta_i -\bar{t}\}$ is a set of discontinuity moments of  $\xi(t)$  and  the function is a solution of  (\ref{eq:graz1}) \cite{ref1}.
    The next step is to show that the following equality $x(-t ,0,x(t ,0,x_0))=x_0,$ holds for all $x_0 \in D$ and $t\in \mathbb{R}.$ Consider the case $t>0.$ If the set of discontinuity moments $\{\theta_i\}$ is empty, the proof is same with that for continuous dynamical systems \cite{Hirsh-smale}.  Because of the condition $(C2),$ which corresponds to invertibility of the jump function $J,$ the equality  $x(\theta_i ,0, x(\theta_i+))=x(\theta_i),$ holds for all $i\in\mathscr{A}.$ Assuming that $\theta_{-1}<0<\theta_1,$ we should verify $x(-\theta_1,0, x(\theta_1, 0,x_0))=x_0.$ Denote by $\bar{x}(t)=x(t,0,x(\theta_1)).$ The point $x(\theta_1)$ lies on the discontinuity surface $\Gamma.$ By condition $(C3)$ the solution $\bar{x}(t)$ is a trajectory of $y'=f(y)$ for decreasing $t.$ Condition $(C10),$ part $(a),$ implies that the trajectory $\bar{x}(t)$ cannot meet with $\tilde{\Gamma}$ if $t>-\theta_1$ as time decreases.   That   is, $\bar{x}(-\theta_1) = x_0$  as the   dynamics is continuous.     The proof for $t<0$ can be done in a similar way. $\square$

\begin{remark}
For the application of the results, it is possible to take the initial moment   as $t_0 =0,$ without being the discontinuity moment since of the group property. Then  $x_0\notin\Gamma\cup\tilde{\Gamma}.$
\end{remark}

Denote by $\widehat{[a,b]},$ $a,b\in \mathbb{R},$ the interval $[a,b],$ whenever $a\leq b$ and $[b,a],$ otherwise. Let  $x_1(t)\in PC(\mathbb{R}_{+},\theta^1),$ $\theta^1=\{\theta_i^1\},$ and $x_2 (t)\in PC(\mathbb{R}_{+},\theta^2 ),$ $\theta^2=\{\theta_i^2\},$ be two different solutions of (\ref{eq:graz1}).

\begin{definition} The solution $x_2(t)$ is in the $\epsilon-$neighborhood of $x_1(t)$ on the interval $\mathscr{I}$ if 
\begin{itemize}
\item the sets $\theta^1$ and $\theta^2$ have same number of elements in $\mathscr{I};$ 
\item $|\theta_i^1-\theta_i^2|<\epsilon$ for all $\theta_i^1\in \mathscr{I};$
\item the inequality $||x_1(t)-x_2(t)||<\epsilon$ is valid for all t, which satisfy $t\in \mathscr{I}\setminus\cup_{\theta_i^1\in \mathscr{I}} (\theta_i^1-\epsilon,\theta_i^1+\epsilon).$ 
\end{itemize}
\end{definition}

The topology defined with the help of $\epsilon-$ neighborhoods is called the B-topology. It can be apparently seen that it is Hausdorff and it can be considered also if two solutions
$x_1(t)$ and $x_2(t)$ are defined on a semi-axis or on the entire real axis.

\begin{definition} The solution $x_0 (t)=x(t,0,x_0),\, t\in\mathbb R, \, x_0\in D,$ of (\ref{eq:graz1}) B-continuously depends on $x_0$ for increasing t if there corresponds  a positive number $\delta$ to any positive $\epsilon$ and a finite interval $[0,b],\ b>0$ such that any other solution $x(t)=x(t,0, \tilde{x} )$ of (\ref{eq:graz1}) lies in $\epsilon-$neighborhood of $x_0 (t)$ on $[0,b]$ whenever $\tilde{x} \in B(x_0,\delta).$ Similarly, the solution $x_0 (t)$ of (\ref{eq:graz1}) B-continuously depends on $x_0$ for decreasing t if there corresponds  a positive number $\delta$ to any positive $\epsilon$ and a finite interval $[a,0],\ a<0$ such that any other solution $x(t)=x(t,0,\tilde{x} )$ of (\ref{eq:graz1}) lies in $\epsilon-$neighborhood of $x_0 (t)$ on $[a,0]$ whenever $\tilde{x} \in B(x_0,\delta).$ The solution $x_0 (t)$ of (\ref{eq:graz1}) B-continuously depends on $x_0$ if it continuously depends on the initial value, $x_0,$ for both increasing and decreasing $t.$
\end{definition}

If conditions (C1)-(C7) hold, then each solution $x_0 (t):\mathbb{R}\rightarrow\mathbb{R}^n,$ $x_0 (t)=x(t,0,x_0),$ of (\ref{eq:graz1}) continuously depends on $x_0$ \cite{ref1}.

\subsection{B-equivalence  to  a system  with  fixed moments of impulses}\label{B-equivalence}

In order to facilitate the analysis of the system with variable moments of impulses \eqref{eq:graz1},  a \textit{B-equivalent system} \cite{ref1} to the system with variable moments of impulses will be utilized in our study. Below, we will construct the B-equivalent system. 

Let $x(t)=x(t,0,x_0+\Delta x)$ be a solution of system (\ref{eq:graz1}) neighbor to $x_0(t)$ with small $\|\Delta x\|.$ If the point $x_0(\theta_i)$ is a $(\beta)-$   or  $(\gamma)-$ type point, then it   is a boundary point.  For this reason, there exist two different possibilities for the near  solution $x(t)$  with respect to the surface of discontinuity. They are:

\begin{itemize} \label{conditions}
\item[$(N1)$] The solution $x(t)$ intersects the surface of discontinuity, $\Gamma,$ at a moment near to $\theta_i,$
\item[$(N2)$] The solution $x(t)$ does not intersect  $\Gamma,$ in a  small time interval centered at $\theta_i.$ 
\end{itemize}



Consider a solution $x_{0}(t): \mathscr{I}\rightarrow \mathbb{R}^n,$ $\mathscr{I}\subseteq \mathbb{R},$ of (\ref{eq:graz1}). Assume that all discontinuity points $\theta_i,$ $i\in\mathscr{A}$ are interior points of $\mathscr{I}.$   There exists a positive number $r,$ such that $r$-neighborhoods of $D_i(r)$ of $(\theta_i ,x_0 (\theta_i))$ do not intersect each other.  Consider $r$ is sufficiently small and so that every solution of \eqref{ode} which   satisfies    condition  $(N1)$  and  starts in $D_i (r)$ intersects $\Gamma$ in $G_i (r)$    as $t$ increases or decreases.
Fix $i\in\mathscr{A}$ and let $\xi(t)=x(t,\theta_i,x),$ $(\theta_i,x)\in D_i (r),$ be a solution of \eqref{ode}, $\tau_i =\tau_i (x)$ the meeting time of $\xi(t)$ with $\Gamma$ and $\psi(t)=x(t,\tau_i,\xi(\tau_i)+J(\xi(\tau_i)))$ another solution  of \eqref{ode}.  Denoting by $W_i(x)=\psi(\theta_i )-x,$ one can find that it is equal to
\begin{equation}\label{W-map}
W_i (x)= \int_{\theta _i}^{\tau _i}{f(\xi(s))ds}+J(x+\int_{\theta _i}^{\tau _i}{f(\xi(s))ds})+\int_{\tau _i}^{\theta _i}{f(\psi(s))ds}
\end{equation}
and maps an intersection of the plane $t=\theta_i$ with $D_i (r)$ into the plane $t=\theta_i.$ 

Let us present the following system of differential equations with impulses at fixed moments, whose impulse moments, $\{\theta_i\},\ i\in\mathscr{A},$ are the moments of discontinuity of $x_0 (t),$

\begin{equation}\label{eq:graz2}
\begin{aligned}
&y'=f(y),\\
&\displaystyle{\Delta y|_{t=\theta_i }= W_i (y(\theta_i))}.
\end{aligned}
\end{equation}

The function $f$ is the same as the function in system (\ref{eq:graz1}) and the maps $W_i,$ $i\in\mathscr{A},$ are defined by equation (\ref{W-map}). If    $\xi(t)=x(t,\theta_i,x)$  does not   intersect 
$\Gamma$  near   $\theta_i$  then we   take  $W_i(x)  = 0.$    

Let us introduce the sets $F_r =\{(t,x)|t\in I, \|x-x_0 (t)\|<r\},$ and $\bar{D}_i (r),$ $i\in\mathscr{A},$   closure of an $r-$ neighborhood of the point $(\theta_i ,x_0 (\theta_i+)).$ Write $D^r = F_r \cup (\cup_{i\in\mathscr{A}} D_i (r))\cup (\cup_{i\in\mathscr{A}} \bar{D}_i (r)).$ Take $r>0$ sufficiently small so that $D^r \subset\mathbb{R}\times D.$ Denote by $D(h)$ an $h$-neighborhood of $x_0 (0).$
 Assume that conditions  $(C1)-(C10)$  hold. Then systems
(\ref {eq:graz1}) and (\ref {eq:graz2}) are B-equivalent in $D^{r}$ for 
 a sufficiently small $r$ \cite{ref1}. That is, if there exists  $h >0,$  such that:
\begin{enumerate}
\item for every solution $y(t)$ of (\ref{eq:graz2}) such that $y(0) \in D(h),$ the
integral curve of $y(t)$ belongs to $D^{r}$ and     there exists
a solution
   $x(t) = x(t,0,y(0))$ of    (\ref{eq:graz1}) which satisfies
\begin{eqnarray}
&&x(t) = y(t), \ \ t \in [a,b]\backslash \cup_{i=-k}^{m} (\widehat
{\tau_i,\theta_i]}, \label{e12}
\end{eqnarray}
where $\tau_i$ are moments of discontinuity of $x(t).$  One   should   precise    that    we   assume  $\tau_i =  \theta_i,$ if $x(t)$ satisfies $(N2).$
Particularly,
\begin{equation}
\begin{array}{l}
x(\theta_i)  = \left \{\begin{array} {ll} y(\theta_i), \quad \mbox {\, if $ \theta_i \leq \tau_i$},\\
y(\theta_i^+), \mbox {otherwise,} \end{array}\right.   \\
y(\tau_i)  = \left \{\begin{array} {ll} x(\tau_i), \quad \mbox {\, if $ \theta_i \geq\tau_i$},\\
x(\tau_i^+), \mbox {otherwise.} \end{array}\right. \label{e13}
\end{array}
\end{equation}
\item Conversely, if  (\ref{eq:graz2}) has a solution
$y(t)=y(t,0,y(0)), y(0) \in D(h),$ then there exists a solution
$x(t)=x(t,0,y(0))$ of (\ref{eq:graz1}) which has an integral curve  in
$D^{r},$ and (\ref{e12})  holds. \label{deqv}
\end{enumerate}
\label{defneq}

A solution $x_0(t)$ satisfies (\ref {eq:graz1}) and (\ref {eq:graz2})
simultaneously.

%
Consider a solution $x_0(t):\mathbb{R}\rightarrow\mathbb{R}^n,$ $x_0(t)=x(t,0,x_0), \ x_0\in D$ with discontinuity moments $\{\theta_i\}$. Fix a discontinuity moment $\theta_i.$  At this discontinuity moment, the trajectory may be on $\Gamma$ and $\tilde{\Gamma}.$ All possibilities of discontinuity moment should be analyzed. For this reason,  we should investigate the following six cases:

\begin{itemize}
\item[$(\alpha)$] $x_0(\theta_i)\in\Gamma\setminus\partial\Gamma$, \quad \quad \quad \quad \quad \quad \quad \quad \quad \quad \quad \quad \quad \quad  $(\alpha^{\prime})$ $x_0(\theta_i)\in\tilde{\Gamma}\setminus\partial{\tilde{\Gamma}}$, 
\item[$(\beta)$] $x_0(\theta_i)\in\partial{\Gamma}$ $\&$  $\langle \nabla\Phi(x_0(\theta_i)),f(x_0(\theta_i))\rangle \neq 0,$  \quad $(\beta^{\prime})$ $x_0(\theta_i)\in\partial{\tilde{\Gamma}}$ $\&$  $\langle \nabla\tilde{\Phi}(x_0(\theta_i)),f(x_0(\theta_i))\rangle\neq0, $
\item[$(\gamma)$] $x_0(\theta_i)\in\partial{\Gamma}$ $\&$  $\langle \nabla\Phi(x_0(\theta_i)),f(x_0(\theta_i))\rangle = 0, $  \quad $(\gamma^{\prime})$ $x_0(\theta_i)\in\partial{\tilde{\Gamma}}$ $\&$  $\langle \nabla\tilde{\Phi}(x_0(\theta_i)),f(x_0(\theta_i))\rangle=0. $
\end{itemize}

If a discontinuity point $x_0(\theta_i)$ satisfy the case $(\alpha),\  ((\alpha^{\prime}))$ the case $(\beta),\ ((\beta^{\prime}))$ and the case $(\gamma), ((\gamma^{\prime}))$  we will call it an $(\alpha)-$ type point, a $(\beta)-$ type point and a $(\gamma)-$ type point, respectively.

Besides, we present the following definition which is compliant with Definition \ref{grazingorbit}.

\begin{definition} \label{grazingsolution}
If there exists  a discontinuity moment, $\theta_i,$ $i\in\mathscr{A},$ for which one of the cases $(\gamma)$ or $(\gamma^{\prime})$ is valid, then the solution $x_0(t)=x(t,0,x_0),$ $x_0\in\mathbb{R}^n$ of (\ref{eq:graz1}) is called \textit{a grazing solution}    and  $t=\theta_i$ is called a \textit{grazing moment}.  
\end{definition}

Next, we consider \textit{the differentiability properties of grazing solutions}. The theory for  the smoothness of   discontinuous dynamical systems' solutions without grazing phenomenon is provided in \cite{ref1}.


Denote by $\bar x(t), j=1,2,\ldots,n,$   a solution  of   (\ref{eq:graz2}) such that  $\bar x(0) = x_0 + \Delta x,\, \Delta x = (\xi_1,\xi_2,\ldots,\xi_n),$  and  let $\eta_i$  be the moments of discontinuity of $\bar x(t).$

The following conditions are required in what follows.

\begin{itemize}
\item[(A)] For all $t \in [0,b] \backslash \cup_{i\in\mathscr{A}} \widehat
{(\eta_i,  \theta_i]}, $ the following  equality is satisfied
\begin{eqnarray}
&& \bar x(t)  - x_0(t) = \sum\limits^{n}_{i=1} u_i(t) \xi_i + O(\|\Delta x\|), \label{edd2}
\end{eqnarray}
where $u_i(t) \in PC([0,b],\theta).$

\item[(B)] There exist constants $\nu_{ij}, j\in\mathscr{A},$ such  that
\begin{eqnarray}
&&  \eta_j - \theta_j = \sum\limits^{n}_{i=1} \nu_{ij} \xi_i + O(\|\Delta x\|); \label{edd1}
\end{eqnarray}

\item[$(C)$] The discontinuity moment $\eta_j$ of the near solution  approaches to the discontinuity moment  $\theta_j, j\in\mathscr{A},$ of grazing one as $\xi$ tends to zero.
\end{itemize}

The   solution $\bar x(t)$   has  a linerization   with  respect   to  solution $x_0(t)$   if   the   condition  $(A)$ is  valid  and, moreover,   if the  point  $x_0(\theta_i)$ is of $(\alpha)-$ or $(\beta)-$ type, then the condition $(B)$ is fulfilled.  For the case $x_0(\theta_i)$ is of $(\gamma)-$ type the condition $(C)$ is true.    
 
 The  solution $x_0(t)$ is $K-$differentiable  with  respect to    the   initial  value $x_0$  on $[0,b]$ if  for  each  solution $\bar x(t)$   with    sufficiently  small  
 $\Delta x$   the   linearization exists. The   functions $u_i(t)$  and   $\nu_{ij}$   depend   on  $\Delta x$   and   uniformly  bounded  on a  neighborhood of  $x_0.$  
 

 It is easy to see that the differentiability implies $B-$continuous dependence on solutions to initial data.

Define the map $\zeta(t,x)$  as $\zeta(t,x)=x(t,0,x),$ for  $x\in D.$

A $K$-smooth discontinuous flow is a map $\zeta(t,x):\mathbb{R}\times D\rightarrow D,$ which satisfies the following properties:
\begin{itemize} 
\item[(I)] The group property:
\begin{itemize}
\item[(i)] $\zeta(0,x):D\rightarrow D$ is the identity;
 \item[(ii)] $\zeta(t,\zeta(s,x))=\zeta(t+s,x)$ is valid for all $t,s\in\mathbb{R}$ and $x\in D.$
\end{itemize}
\item[(II)] $\zeta(t,x)\in PC^1 (\mathbb{R})$ for each fixed $x\in D.$
\item[(III)] $\zeta(t,x)$ is $K$-differentiable in $x\in D$ on $[a,b]\subset\mathbb{R}$ for each $a,b$ such that the discontinuity points of $\zeta(t,x)$ are interior points of  $[a,b].$
\end{itemize}

In \cite{ref1}, it was proved that if  the conditions of Theorem \ref{continue1}  and (C1)-(C10)  are fulfilled, then system (\ref{eq:graz1}) defines a $B$-smooth discontinuous flow \cite{ref1}   if there   is no    grazing points
for   the   dynamics.   It is easy to observe   that the   $B$-smooth discontinuous flow  is a subcase  of the  $K$-smooth discontinuous flow.
In the next section, we will construct a variational system for (\ref{eq:graz1}) in the neighborhood of grazing orbits. That is, we will assume that some of the discontinuity points are $(\gamma)-$ type   points.  Linearization around a solution and its stability will be taken into account. Thus,  analysis  of the   discontinuous   dynamical    systems   with   grazing points will  be completed. 

\section{Linearization around grazing  orbits  and discontinuous dynamics}

The object of this section is to verify $K -$ differentiability of the grazing solution. Consider a grazing solution $x_0(t)=x(t,0,x_0),\ x_0\in D,$ of (\ref{eq:graz1}). 
We will demonstrate that one can write the variational system for the solution $x_0(t)$ as follows:

\begin{equation}\label{lin1}
\begin{aligned}
&u'=A(t)u, \\
&\Delta u|_{t=\theta_i} =B_i u(\theta_i),
\end{aligned} 
\end{equation} 

where the matrix $A(t)\in \mathbb{R}^{n\times n}$ of the form $A(t)=\frac{\partial f(x_0(t))}{\partial x}.$   The matrices $B_i, \ i=1,\ldots,n,$   will be defined in the remaining part of the paper. The matrix $B_i$ is bivalued if $\theta_i$ is a grazing moment  or   of $(\beta)-$type.  

The right hand side of the second  equation in (\ref{lin1}) will be described in the remaining part of the paper for each type of the points. As the \textit{linearization at a point of discontinuity}, we comprehend the second  equation in (\ref{lin1}).

\subsection{Linearization at $(\alpha)-$type points} \label{apoint}

Discontinuity   points of   $(\alpha)$ and  $(\alpha^{\prime})$  types   are   discussed in \cite{ref1}.  In this subsection, we will   outline    the   results of the book.

 Assume that $x(\theta_i)$ is an $(\alpha)-$type point. It is clear that the $B-$ equivalent system (\ref{eq:graz2}) can be applied in the analysis. The functions $\tau_i(x)$ and  $W_i(x),$ are described in Subsection \ref{B-equivalence}. Differentiating $\Phi(x(\tau_i(x)))=0,$ we have

\begin{equation}\label{dert}
\frac{\partial\tau_i(x_{0}(\theta_i))}{\partial x_j}=-\frac{\Phi_{x}(x_0(\theta_i))\frac{\partial x_{0}(\theta_i)}{\partial x_{0j}}}{\Phi_{x}(x_0(\theta_i))f(x_0(\theta_i))}.
\end{equation}

Then, considering (\ref{W-map}), we get the following equation,
\begin{equation}\label{derw}
\frac{\partial W_{i}(x_0(\theta_i))}{\partial x_{0j}}=(f(x_0(\theta_i))-f(x_0(\theta_i)+J(x_0(\theta_i))))\frac{\partial\tau_i}{\partial x_{0j}}+\frac{\partial I}{\partial x}(e_j+f\frac{\partial\tau_i}{\partial x_{0j}}),
\end{equation}
where $e_j=(\underbrace{0,\ldots,1}_j,\ldots,0).$ 

The matrix $B_i\in \mathbb{R}^{n\times n}$ in equation (\ref{lin1}) is defined as $B_i=W_{ix},$ where  $W_{ix}$ is the $n\times n$ matrix of the form $W_{ix}=[\frac{\partial W_{i}(x_0(\theta_i))}{\partial x_1},\frac{\partial W_{i}(x_0(\theta_i))}{\partial x_2},\ldots,\frac{\partial W_{i}(x_0(\theta_i))}{\partial x_n}].$ Its vector-components $\frac{\partial W_{i}(x_0(\theta_i))}{\partial x_{0j}}, \quad j=1,\ldots,n,$  evaluated by (\ref{derw}). Moreover, the components of the gradient $\nabla \tau_i$ have to be evaluated by formula (\ref{dert}).  

\subsection{Linearization at $(\beta)-$type points} \label{bpoint}

In what follows, denote $n\times n$ zero matrix by $O_{n}.$ In the light of the possibilities $(N1)$ and $(N2),$ the matrix $B_i $ in (\ref{eq:graz1}) can be expressed as follows:

\begin{eqnarray} \label{beta} B_i=\begin{cases}  O_{n}, \quad &\mbox{if } \quad \mbox{ $(N1)$ is valid,} \\ 
 W_{ix}, \quad &\mbox{if } \quad \mbox{ $(N2)$ is valid,}  \end{cases}
\end{eqnarray}
where $ W_{ix}$  is  evaluated by formula (\ref{derw}) and $\nabla\tau(x)$ evaluated by formula \eqref{dert}. 

The differentiability properties for the cases $(\alpha^{\prime})$ and $(\beta^{\prime})$ can be investigated similarly.


\subsection{Linearization at a grazing point} \label{graz1}

Fix a discontinuity  moment  $\theta_i$ and assume that one of the cases $(\gamma)$ or $({\gamma^{\prime}})$ is satisfied. We will investigate the case $(\gamma).$ The case $({\gamma^{\prime}})$ can be considered in a similar way. 

Considering condition $(C1)$ with the formula  \eqref{dert}, it is easy to see that one coordinate of it is infinity at a grazing point. This gives arise singularity in the system, which makes the analysis harder  and  the  dynamics   complex.  Through the formula \eqref{dert}, one can see that the singularity is just caused by the position of the vector field  with respect to the surface of discontinuity  and the impact component of the dynamical system does not participate in the appearance of the singularity.  To handle with the singularity, we will rely on the following conditions.

\begin{itemize}
\item[$(A1)$]  A  grazing   point  is isolated.   That  is, there  is a neighborhood  of  the  point with  no  other   grazing points.  

\item[$(A2)$] The   map  $W_i(x)$  in   \eqref{W-map} is   differentiable  at the grazing point $x = x_0(\theta_i).$

\item[$(A3)$]  The function $\tau_i(x)$  does not exceed   a positive  number  less than  $\theta_{i+1} -\theta_i$    near a grazing point,  $x_0(\theta_i),$ on a set  of points     which   satisfy   condition $(N1).$

\end{itemize}


In the present  paper,  we analyze   the  case,   when   the impact  functions neutralize  the  singularity caused by transversality. That is,  the triad: impact law, the surface of discontinuity and the vector field is specially chosen,   such    that  condition $(A2)$  is valid.  Presumably,   if there   is no of this  type  of  suppressing,  complex   dynamics near the  grazing  motions may  appear \cite{budd2001,Luo2006,Nord1,Nord97}.
In the examples stated in the remaining part of the paper, one can see the verification of $(A2),$ in details.

Let   us prove the   following assertion.   

\begin{lemma}\label{conttau}   If   conditions $(C1),$ $(C4),$ $(C6),$ $(C8)$ and $(A3)$ hold. Then, $\tau_i(x)$  is continuous near   a grazing point $x_0(\theta_i),$ on a set  of points,   which   satisfy   condition $(N1).$  
\end{lemma}

\textbf{ Proof.} Let $x_0(\theta_i)$ be a grazing point. If $\bar x$ is not a point from the orbit of the grazing solution, the continuity of $\tau_i(x)$ at the point $x=\bar x$ can be proven using similar technique presented in \cite{ref1}. Now, the continuity at $x_0(\theta_i)$ is taken into account. On the contrary, assume that $\tau_i(x)$ is not continuous at the point $x=x_0(\theta_i).$ Then, there exists a positive number $\epsilon_0$ and a sequence $\{x_n\}_{n\in\mathbb{Z}}$ such that $\tau_i(x_n)>\epsilon_0$ whenever $x_n\rightarrow x_0(\theta_i),$ as $n\rightarrow \infty.$ Moreover, from condition $(A3),$  one can assert that there exists a subsequence $\tau_i(x_{n_k})$ which converges to a number $\epsilon_0\leq \tau_0  <  \theta_{i+1} -\theta_i.$ Without loss of generality, assume that the subsequence converges the point where the sequence $\{x_n\}_{n\in\mathbb{Z}}$ converges.  Since of the continuity of solutions in initial value,  $x(\tau_i(x_n),0,x_n)$ approaches to $x(\tau_0,0,x(\theta_i)).$ But $x(\tau_i(x_n),0,x_n)$ is on the surface of discontinuity $\Gamma,$ $x(\tau_0,0,x_0(\theta_i))\notin \Gamma.$ This contradicts with the closeness of the surface of discontinuity $\Gamma.$  The continuity at other points of the grazing orbit is valid by the group property. $\square$

Since   of $B-$equivalence  of systems (\ref{eq:graz1}) and  (\ref{eq:graz2}),  we   will   consider   linearization    around $x_0(t)$  as solution of the   system (\ref{eq:graz2}),  consequently, only   formula   \eqref{edd2}   will be needed.    Finally, the linearization matrix   for the   grazing point  also has to  be defined by the formula \eqref{beta}, where $W_{ix}$ exists by condition $(A2).$

In what follows, we will consider only grazing motions such that   condition $(A2)$  holds.  Consequently,  the continuous dependence  on initial data is valid. More precisely,  $B-$ continuous dependence on initial data is true.  Now, if conditions $(C1)-(C10)$   and $(A1),(A2)$  are assumed,  the system \eqref{eq:graz1} defines a $K-$ smooth discontinuous flow for dynamics with grazing points.

\subsection{Linearization around a grazing periodic solution}
Let $\Psi(t):\mathbb{R}\rightarrow D$ be a periodic solution of (\ref{eq:graz1})  with   period $\omega>0$ and $\theta_i,\ i\in\mathbb{Z},$ are the points of discontinuity which satisfy  $(\omega,p)-$ property, i.e. $\theta_{i+p}=\theta_i + \omega, p$  is a natural number.   

%
%
%
%
%
%

Let us fix a solution $x(t)=x(t,0,\Psi(0)+\Delta x)$ and assume that linearization of $\Psi(t)$  with   respect   to  $x(t)$  exists and is of the form
\begin{equation}\label{eq:graz1lin1}
\begin{aligned}
&u'=A(t)u,\\
&\Delta u|_{t=\theta_i}= B_{i} u.
\end{aligned}
\end{equation}
The matrix $B_{i}$   is determined by  \eqref{beta}.
It is  known   that   $ A(t+\omega)=A(t), \, t  \in \mathbb R.$   But, the sequence $B_i$ may not be periodic in general, since of \eqref{beta}. This makes the analysis of the neighborhood of $\Psi(t)$ difficult. For this reason, we suggest the following condition. 
 
\begin{itemize}
\item[(A4)] For each   sufficiently   small $\Delta x \in \mathbb{R}^n,$ the variational system  \eqref{eq:graz1lin1}  satisfies $B_{i+p}=B_{i} ,$ $i\in\mathbb{Z}.$ 
There   exist a finite number $m \le 2^l,$ where  $l$ is the number of points of $(\beta)-$ or $(\gamma)-$ type in the interval $[0,\omega],$ of the periodic sequences  $B_{i}.$       
\end{itemize}

The assumption $(A4)$ is valid for many  low   dimensional  models   of  mechanics and  those  which  can be decomposed into  low dimensional subsystems.  To distinguish periodic sequences $B_i$  in the assumption $(A4),$   we will apply the notation $B_i=D_i^{(j)},$ $i\in\mathbb{Z}$ and $j=1,2,\ldots,m.$

If the condition $(A4)$ is not  fulfilled, then   complex dynamics near   a  periodic   motion may  appear.  This   case  can   be investigated   either   by   methods   developed   through  mappings  applications  \cite{nusse94, Piir-virgin-Champhneys} or  it  requests   additional   development of our   present  results.

%

In the next example, we will demonstrate that the system constitutes  $K-$ smooth discontinuous flow although it has grazing points in the phase space. 
\begin{example} \label{exdds}(K-smooth discontinuous flow with grazing points). Consider  an   impact model
\begin{subequations}\label{orbex2}
\begin{equation}\label{orbex2:1}
\begin{aligned}
&y_1 '=y_2 , \\
&y_2 '=-y_1+0.001y_2, 
\end{aligned}
\end{equation}
\begin{equation}\label{orbex2:2}
\begin{aligned}
&\displaystyle{\Delta y_2|_{y\in \Gamma_1}}=-y_2 -R_1y_2^2,\\
&\displaystyle{\Delta y_2|_{y\in \Gamma_2}}= -(1+R_2)y_2,
\end{aligned}
\end{equation}
\end{subequations}
with the domain   $D=\mathbb{R}^2,$ $R_1=\exp(-0.0005\pi)$ and $R_2=0.9.$ In the paper \cite{repeated}, it is stated that the coefficient of restitution for low velocity impact still remains as an open problem. In the study \cite{impactdef}, by considering Kelvin-Voigt model for the elastic impact, we derived quadratic terms of the velocity in the impact law.  This arguments   make the quadratic term for the impulse equation   \eqref{orbex2:2}  reasonable.

Let us describe the set of discontinuity curves by $\Gamma=\Gamma_1\cup\Gamma_2.$ The components $\Gamma_1$ and $\Gamma_2$ are intervals of the vertical lines  $y_1=\exp(0.00025\pi)$ and $y_1=0,$  respectively  and they   will be precised  next. Fix a point $P=(0,\bar{y}_2)\in D,$ with $\bar{y}_2>1.$  Let $y(t,0,P)$ be a solution of (\ref{orbex2:1}) and it  meets with the vertical line $x_1=\exp(0.00025\pi), \ x_2>0$  at the point $P_2=(\exp(0.00025\pi),y_2(\theta_1,0,P)),$ where $\theta_1$ is the meeting moment with the line. Consider the point $Q_2=(\exp(0.00025\pi),-R_1y_2(\theta_1,0,P_2)^2) $ and    denote  $Q_1=(0,y_2(\theta_2,0,Q_2)),$ where $\theta_2$ is the moment of meeting of the solution $y(t,0,Q_2)$ with the vertical line $x_1=0, \ x_2<0.$ We shall need also the point $P_1=(0,-R_2y_2(\theta_2,0,Q_2)).$ Finally,  we   obtain the region $G$ in yellow and blue between   the vertical lines and graphs  of the solutions in Figure \ref{R-region}. The region $G$ contains discontinuous trajectories and outside of this region all trajectories are continuous. Moreover, both region $G$ and its complement are invariant.

\begin{figure}[ht] 
\centering
\includegraphics[width=9 cm]{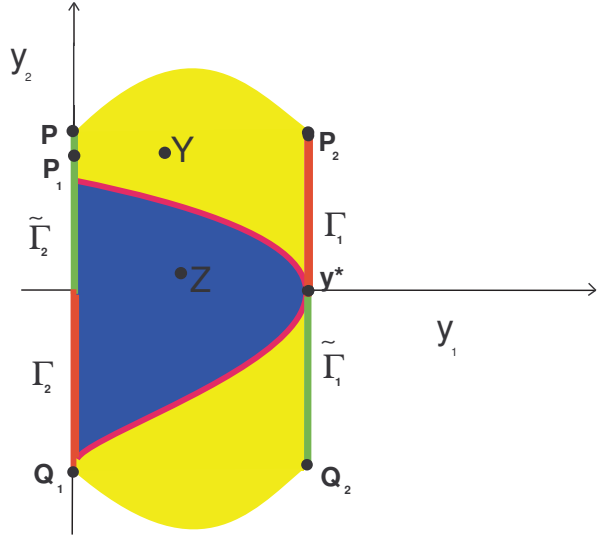}
\caption{The   region  $G$  for  system \eqref{orbex2} is depicted  in details. The curves of discontinuity $\Gamma=\Gamma_1\cup\Gamma_2$  and $\tilde{\Gamma}=\tilde{\Gamma}_1\cup\tilde{\Gamma}_2$  are drawn as vertical lines in red and green, respectively and the grazing orbit  in magenta.}
\label{R-region}
\end{figure}

Define    $ \Gamma_1 =\{(y_1, y_2)| \ y_1=\exp(0.00025\pi), \quad 0\leq y_2\leq y_2(\theta_1,0,(0,\bar{y}_2))\},$  and $\Gamma_2 =\{(y_1, y_2)|\ y_1=0,\, y_2(\theta_2,0,-R_1y_2((\theta_1,0,(0,\bar{y}_2)))^2)\leq y_2\leq 0\}.$    The boundary of the curve, $\Gamma=\Gamma_1\cup\Gamma_2,$ has   of four points, they are $$\partial\Gamma=\{(0,0), (\exp(0.00025\pi),0),(\exp(0.00025\pi),y_2(\theta_1,0,(0,\bar{y}_2))), (0,y_2(\theta_2,0,-R_1y_2(\theta_1,0,(0,\bar{y}_2)))\}.$$  In the following part of the example,  we will show that two of them, $y^*=(y_1^*,y_2^*)=(\exp(0.00025\pi),0)$ and the origin, $(0,0)$ are grazing points. Moreover, it can be easily validated that other   two   points are of $\beta-$type. 

 Issuing from system (\ref{orbex2}), the curve of discontinuity $\tilde{\Gamma}$ consists of two components $\tilde{\Gamma}_1$ and $\tilde{\Gamma}_2.$ The components  are the following sets $$\tilde{\Gamma}_1 =\{(y_1, y_2)| \ y_1=\exp(0.00025\pi), \, -R_1y_2(\theta_1,0,(0,\bar{y}_2))^2\leq y_2\leq 0\}$$ and  $$\tilde{\Gamma}_2 =\{(y_1, y_2)|\ y_1=0,\, 0\leq -R_2y_2(\theta_2,0,Q_2)\}.$$

One can verify that the function
\begin{eqnarray}\label{per}
 \Psi(t) =
  \begin{cases}
  \exp(0.0005 t) \Big(\sin(t), \cos(t)\Big), & \text{if } t \in [0,\pi), \\
    
			(0,1) , & \text{if } t=\pi,\end{cases}
\end{eqnarray}
is a discontinuous periodic solution of (\ref{orbex2}) with period $\omega=\pi,$ whose discontinuity points $(0,1)$ and $(0,-\exp(0.0005\pi))$ belong to $\tilde{\Gamma}$ and $\Gamma,$ respectively. The expression $$\langle \nabla\Phi((\exp(0.00025\pi),0)),f((\exp(0.00025\pi),0))\rangle=\langle(1,0),(0, -\exp(-0.00025\pi))\rangle=0$$ verifies that $y^*$ is a $(\gamma)-$ type point, i.e. a grazing point of the solution $\Psi(t).$  It  is easily  seen that   the grazing   is axial. Now, we can assert that the periodic solution (\ref{per}) is a grazing solution in the sense of Definition \ref{grazingsolution}. Its simulation is depicted in Figure \ref{periodicdisc}.
\begin{figure}[ht] 
\centering
\includegraphics[width=9 cm]{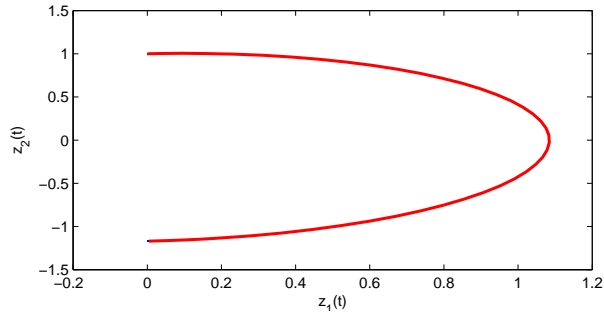}
\caption{The   grazing orbit of   system \eqref{orbex2}.}
\label{periodicdisc}
\end{figure} 

Since   the complement   of $G$  is  invariant  in both   directions and consists  of continuous   trajectories  of   the   linear   system  (\ref{orbex2:1}), one can  easily conclude     that   the   complement  is    a continuous   dynamical   system \cite{Hirsh-smale}.  Thus,   to    verify   the   dynamics for the whole system,  one need  to analyze   it  in the region $G.$
  This  set  is bounded, consequently for   solutions  in it conditions  $(C8)$ and $(C9)$  are fulfilled and   by Theorem \ref{continue1}, they admit   $B-$sequences    and continuation  property.

Consider   a  function  $\zeta(y_2):[y_2(\theta_2,0,Q_2),y_2(\theta_1,0,P)]\rightarrow [y_2(\theta_2,0,Q_2),y_2(\theta_1,0,P)]$  such that   it is continuously differentiable,  satisfies  $\zeta(y_2)=-R_2y_2$  in a neighborhood of $y_2=0$ and is the identity at the boundary points, i.e. $\zeta(y_2(\theta_1,0,P))=y_2(\theta_1,0,P)$ and  $\zeta(y_2(\theta_1,0,P))=y_2(\theta_1,0,P).$ 
It   is easily   seen that such function exists.  On the basis of this discussion, let us introduce   the  following system, 

\begin{equation}\label{orbex22}
\begin{aligned}
&y_1 '=y_2 , \\
&y_2 '=-y_1+0.001y_2, \\
&\displaystyle{\Delta y_2|_{y\in \Gamma}}=\zeta(y_2)-y_2.
\end{aligned}
\end{equation}

It is apparent that   system (\ref{orbex22})  is equivalent to (\ref{orbex2}) near the orbit of  periodic solution $\Psi(t).$   That   is,   they  have the same trajectories there.  

Specifying    (\ref{eq:graz1})  for (\ref{orbex22}), it is easy to obtain that  $\Phi(y_1,y_2)=\tilde{\Phi}(y_1,y_2)=(y_1-\exp(0.00025\pi))y_1,$   $f(y_1,y_2)=(y_2,-y_1+0.001y_2)$ and  $ J(y) = (y_1,\zeta(y_2)).$

Now, we will verify that system (\ref{orbex22}) defines a $K-$ smooth discontinuous flow. First, condition $(C1)$  is verified since $\nabla\Phi_1(y)=\nabla\Phi_2(y)=(1,0)\neq 0,$ for all $y\in D.$ The jump function $J(y)=(y_1,\zeta(y_2))$ is  continuously differentiable function. So, condition $(C2)$ is valid. It is true that  $\Gamma\cap\tilde{\Gamma}\subseteq\partial\Gamma\cap\tilde{\partial\Gamma}.$    Inequalities  $\langle \nabla\Phi_1(y), f(y)\rangle=\langle (1,0), (y_2, -y_1+0.001y_2)\rangle=y_2\neq 0$  and $\langle \nabla\Phi_2(y), f(y)\rangle=\langle (1,0),  (y_2, -y_1+0.001y_2 )\rangle=y_2\neq 0,$ if $y\in\Gamma\setminus\partial\Gamma,$  validate the condition $(C4).$ Moreover, $\langle \nabla\tilde{\Phi}_1(y), f(y)\rangle=\langle (1,0), (y_2, -y_1+0.001y_2  )\rangle=y_2\neq 0$  and $\langle \nabla\tilde{\Phi}_2(y), f(y)\rangle=\langle (1,0), (y_2, -y_1+0.001y_2)\rangle=y_2\neq 0,$ if $y\in\tilde{\Gamma}\setminus\partial\tilde{\Gamma}.$  Conditions $(C6)$ and $(C7)$  hold as the function $\zeta$ is such defined. 
Thus, conditions $(C1)-(C10)$   have been  verified. Consequently, the system \eqref{orbex2} defines the $K-$ smooth discontinuous flow for all motions except the grazing ones. To complete the discussion, one need to linearize the system near the grazing solutions. First, we proceed with the linearization around the grazing periodic orbit (\ref{per}).

The   solution, $\Psi(t)$  has two  discontinuity moments $\theta_1=\frac{\pi}{2}$ and $\theta_2=\omega$ in the interval  $[0,\omega].$ The   corresponding discontinuity points are of $(\gamma)-$ and $(\alpha)-$ types, respectively. Next, we will linearize the system at these points.   The linearization at the  second  point exists  \cite{ref1}  and  the details of this  will be analyzed in the next example. This time, we will  focus   on the grazing point $y^*.$  


First, we assume that  $y(t)=y(t,0,y^*+\Delta y),$ $\Delta y=(\Delta y_1,\Delta y_2)$ is not a grazing solution. Moreover, the solution intersects the line $\Gamma_1$ at time $t=\xi$ near $t=\theta_1$ as time increases. The meeting   point $\bar y=(\bar y_1,\bar y_2) = (y_1(\xi,0,(y^*+\Delta y)),y_2(\xi,0,(y^*+\Delta y)),$  is transversal one. 
It  is clear $\bar y_1=\exp(0.00025\pi)$ and $\bar y_2>0.$ In order to find a  linearization at the moment $t=\theta_i,$ we use formula \eqref{W-map} for $y(t),$   and find that  
\begin{eqnarray}\label{derW_1}
\displaystyle{\frac{\partial W_i(y)}{\partial y_1^0}}&=\displaystyle{\int\limits_{\theta_i}^{\tau(y)}{\frac{\partial f(y(s))}{\partial y}\frac{\partial y(s)}{\partial y_1^0}ds}}+f(y(s))\frac{\partial\tau(y)}{\partial y_1^0}+J_y(y)(e_1+f(y(s))\frac{\partial\tau(y)}{\partial y_1^0})+f(y(s)+J(y(s)))\frac{\partial\tau(y)}{\partial y_1^0}\nonumber   \\
&\displaystyle{+\int\limits_{\tau(y)}^{\theta_i}{\frac{\partial f(y(s)+J(y(s)))}{\partial x}\frac{\partial y(s)}{\partial y_1^0}ds}},
\end{eqnarray}
where $e_1=(1,0)^T,$ $T$ denotes the transpose of a matrix.  
Substituting $y= \bar y$ to the formula \eqref{derW_1}, we obtain that 
\begin{eqnarray}\label{derW_2}
&\displaystyle{\frac{\partial W_i(y(\xi,0,y^*+\Delta y))}{\partial y_1^0}=f(y(\xi,0,y^*+\Delta y))\frac{\partial\tau(y(\xi,0,y^*+\Delta y))}{\partial y_1^0}}\nonumber\\
&\displaystyle{+J_y(y(\xi,0,y^*+\Delta y))\Bigg(e_1+f( y(\xi,0,y^*+\Delta y)))\frac{\partial\tau(y (\xi,0,y^*+\Delta y)))}{\partial y_1^0}\Bigg)}\nonumber \\
&+\displaystyle{f(y(\xi,0,(J(y(\xi,0,y^*+\Delta y)))))\frac{\partial\tau(J(y(\xi,0,y^*+\Delta y)))}{\partial y_1^0}}. 
\end{eqnarray}

Considering the formula \eqref{dert} for the transversal point $\bar y=(\bar y_1,\bar y_2),$ the first component $\displaystyle{\frac{\partial \tau(\bar y)}{\partial y_1^0}}$ can be evaluated as $\displaystyle{\frac{\partial \tau(\bar y)}{\partial y_1^0}=-\frac{1}{\bar y_2}}.$ From the last equality, it  is seen how the singularity appears at the grazing point.  Finally, we obtain that

\begin{eqnarray}\label{derW_3}
\displaystyle{\frac{\partial W_i(\bar y)}{\partial y_1^0}}= \begin{bmatrix} \bar y_2 \\ -\bar y_1-0.001\bar y_2\end{bmatrix}\displaystyle{\Big(-\frac{1}{\bar y_2}\Big)}+\begin{bmatrix} 1 &0 \\ 0 & -2R_1\bar y_2 \end{bmatrix} \Bigg(e_1+\begin{bmatrix} \bar y_2 \\ -\bar y_1-0.001\bar y_2\end{bmatrix}\displaystyle{\Big(-\frac{1}{\bar y_2}\Big)}\Bigg)
\end{eqnarray}
$$-\begin{bmatrix} -R_1(\bar y_2)^2 \\ -\bar y_1+0.001R_1(\bar y_2)^2 \end{bmatrix}\displaystyle{\Big( -\frac{1}{\bar y_2} \Big)}= \begin{bmatrix} \bar y_2 -R_1(\bar y_2)^2 \\ -\bar y_1-0.001(\bar y_2 -R_1(\bar y_2)^2) \end{bmatrix}\displaystyle{\Big(-\frac{1}{\bar y_2}\Big)}+\begin{bmatrix} 1 &0 \\ 0 & -2R_1\bar y_2 \end{bmatrix}\begin{bmatrix} 0\\ \displaystyle{\frac{\bar y_1+0.001\bar y_2}{\bar y_2}} \end{bmatrix}.$$

Calculating the righthand side of  \eqref{derW_3} we have
\begin{eqnarray}\label{derW_4}
\displaystyle{\frac{\partial W_i(\bar y)}{\partial y_1^0}}= \begin{bmatrix} -R_1\bar y_2-1 \\ 0.001(1-R_1\bar y_2)+2R_1(0.001\bar y_2-\bar y_1)\end{bmatrix}.
\end{eqnarray}

The   last  expression   demonstrates   that  the   derivative  is   a  continuous  function of its arguments   in a neighborhood of the grazing point. Since it is defined and continuous   for the points, which  are not from the grazing  orbit  by  the  last  expression and for   other points it  can be determined by  the limit  procedure.   Indeed, one can easily   show that    the  derivative at  the grazing point  $y^*$ is
  \begin{equation}\label{firstvec}
 \begin{bmatrix} -1 \\ 0.001-1.8\exp(0.00025\pi)\end{bmatrix}.  
  \end{equation}
  
  Similarly,  all other points of the   grazing orbit  can be discussed.

Next, differentiating  \eqref{W-map}  with  $y(t)$    again  we obtain that    
\begin{eqnarray}\label{derW_5}
\displaystyle{\frac{\partial W_i(y)}{\partial y_2^0}}&=\displaystyle{\int\limits_{\theta_i}^{\tau(y)}{\frac{\partial f(y)}{\partial y}\frac{\partial y(s)}{\partial y_2^0}ds}}+f(y(s))\frac{\partial\tau(y)}{\partial y_2^0}+J_y(y)(e_2+f(y(s))\frac{\partial\tau(y)}{\partial y_2^0})+f(y +J(y ))\frac{\partial\tau(y)}{\partial y_2^0}\nonumber   \\
&\displaystyle{+\int\limits_{\tau(y)}^{\theta_i}{\frac{\partial f(y(s)+J(y(s)))}{\partial x}\frac{\partial y(s)}{\partial y_2^0}ds}},\end{eqnarray}
where $e_2=(0,1)^T.$ Calculate the  right hand side  of  \eqref{derW_5} at the point $\bar y=(\bar y_1,\bar y_2)$   to   obtain 
\begin{eqnarray}\label{derW_6}
&\displaystyle{\frac{\partial W_i(y(\xi,0,y^*+\Delta y))}{\partial y_2^0}=f(y(\xi,0,y^*+\Delta y))\frac{\partial\tau(y(\xi,0,y^*+\Delta y))}{\partial y_2^0}}\nonumber\\
&\displaystyle{+J_y(y(\xi,0,y^*+\Delta y))\Bigg(e_2+f(y(\xi,0,y^*+\Delta y)))\frac{\partial\tau(y(\xi,0,y^*+\Delta y))}{\partial y_2^0}\Bigg)}\nonumber \\
&+\displaystyle{f(y(\xi,0,y^*+\Delta y))\frac{\partial\tau(y(\xi,0,y^*+\Delta y))}{\partial y_2^0}}. \end{eqnarray}

To calculate the fraction $\displaystyle{\frac{\partial\tau(y(\xi,0,y^*+\Delta y))}{\partial y_2^0}}$ in \eqref{derW_6}, we apply  formula \eqref{dert} for the transversal point $\bar y=(\bar y_1,\bar y_2).$ The second component $\displaystyle{\frac{\partial \tau(\bar y)}{\partial y_2^0}}$ takes the form $\displaystyle{\frac{\partial \tau(\bar y)}{\partial y_2^0}=0.}$ This and  formula \eqref{derW_6}  imply  

\begin{eqnarray}\label{derW_7}
\displaystyle{\frac{\partial W_i(\bar y)}{\partial y_2^0}}= \begin{bmatrix} 0 \\ -2R\bar y_2\end{bmatrix}.
\end{eqnarray}
Similar to \eqref{firstvec},  one can obtain that

\begin{eqnarray}\label{secondvec}
\displaystyle{\frac{\partial W_i(y^*)}{\partial y_2^0}}= \begin{bmatrix} 0 \\ 0\end{bmatrix}.
\end{eqnarray}

Joining \eqref{firstvec} and \eqref{secondvec}, it can be obtained that
\begin{eqnarray}\label{derW_8}
\displaystyle{ W_{iy}(y^*)}= \begin{bmatrix}  -1  & 0 \\ 0.001-1.8\exp(0.00025\pi)&0\end{bmatrix}.
\end{eqnarray}
The     continuity    of the   derivatives   in a neighborhood of  $y^*$  implies   that  the function   $W$   is differentiable at  the  grazing  point  $y = y^*,$ and the  condition $(A2)$   is valid.

Now, on the basis of the discussion made above, one can obtain the  bivalued  matrix   of coefficients for the grazing point as

\begin{eqnarray*}
 B_1=\begin{cases} O_{2}, \ & \mbox{if (N1) is valid}, \\ 
\begin{bmatrix} -1 & 0 \\ 0.001-1.8\exp(0.00025\pi)& 0\end{bmatrix},\ & \mbox{if (N2) is valid}. \end{cases}
\end{eqnarray*}

The matrix $D_1^{(1)} =O_{2}$ is for  near solutions of (\ref{per}) which are  in the region where $Z$ in, see Fig. \ref{R-region}, and do not intersect the curve of discontinuity $\Gamma_1.$ The matrix 
$$D_1^{(2)}=\begin{bmatrix} -1 & 0 \\  0.001-1.8\exp(0.00025\pi) &0\end{bmatrix}$$
is for  near solutions of (\ref{per}), which  intersects the curve of discontinuity $\Gamma_1.$  They   start  in the subregion, where the point  $Y$  is placed. 
Thus,    the linearization  for  $\Psi(t)$   at   the   grazing  point  exists.   Moreover,   since   another  point  of discontinuity  $(0,\exp(0.0005\pi))$ is not   grazing, the linearization at  the   point   exist  as   well as   linearization  at   points   of continuity \cite{ref1,perko}.   Consequently, there   exist linearization  around  $\Psi(t).$ 
 
 To verify  condition $(A3),$ consider a near solution $y(t)=y(t,0,\bar y)$   to $\Psi(t),$  where $\bar y=(0,\bar y_2), \, \bar y_2>\Psi_2(0) = 1,$ which satisfy the condition $(N1).$  It is true that $\displaystyle{\theta_{i+1}-\theta_i=\frac{\pi}{2} = \frac{\omega}{2}.}$   The   first    coordinate   of  the near solution is $\displaystyle{y_1(t)=\bar y\exp(0.0005t)\sin(t)}$  and \\  $y_1(\frac{\omega}{2}) = \displaystyle{y_1(\frac{\pi}{2})=\bar  y\exp(0.00025\pi)} >\exp(0.00025\pi) =\Psi_1(\frac{\omega}{2}).$ Thus,   the meeting moment of near solution $y(t)$ with the surface of discontinuity is less than  $\frac{\omega}{2}.$ So, it implies that $0< \tau(y)< \frac{\pi}{2}-\epsilon$  for a small number $\epsilon$ if the first   coordinate of $\bar y$   is close to  $\exp(0.00025\pi).$   This  validates  condition $(A3).$ Now,  Lemma \ref{conttau} proves the condition $(C).$

Now, let   us  consider    the   point $(0,0).$    We   have  that $\langle \nabla\Phi((0,0)),f((0,0))\rangle=\langle(1,0),(0,0)\rangle=0.$ That  is the   origin is   a   grazing point.   In the same time it is a fixed point  of the system.   For  this particular grazing  point,   we   can   find the   linearization  directly.  Indeed, all the  near   solutions   satisfy   the linear  impulsive system,
\begin{equation}\label{orbexzero}
\begin{aligned}
&x_1 '=x_2 , \\
&x_2 '=-x_1+0.001x_2, \\
&\displaystyle{\Delta x_2|_{x_1=0}=-(1 + R_2)x_2.}
\end{aligned}
\end{equation}

%
%
%

Consider a solution $x(t)=x(t,0,x_0),$ where $x_0=(x_1^0,x_2^0)\not = (0,0)$ with moments of discontinuity $\theta_i, i \in \mathbb Z,$ then the linearization system for the equation around the equilibrium is

\begin{equation}\label{orbexzerolin}
\begin{aligned}
&u_1 '=u_2 , \\
&u_2 '=-u_1+0.001u_2, \\
&\displaystyle{\Delta u_2|_{t=\theta_i}=-(1 + R_2)u_2.}
\end{aligned}
\end{equation}
Indeed, if $u_1(t),$ $u_1(0)=e_1,$ $u_2(t),$ $u_2(0)=e_2,$ are solutions of \eqref{orbexzerolin}, then one can see that $x(t)-(0,0)=x_1^0u_1(t)+x_2^0u_2(t),$ for all $t\in\mathbb R.$ 


We   have obtained that  linearization   exists for   both    grazing solutions  $\Psi(t),$   and the  equilibrium at  the origin.   Moreover,   conditions $(C1)-(C10)$ are   valid and   all other   solutions  are   B-differentiable  in parameters \cite{ref1}. Thus,  the system \eqref{orbex2} defines a $K-$ smooth discontinuous flow in  the plane. \end{example}

In the next example, we will finalize the linearization   around the grazing solution $\Psi(t).$

\begin{example} \label{linear}(Linearization around the  grazing    discontinuous  cycle).  We continue analysis of   the  last  example,   and   complete  the variational  system  for $\Psi(t).$  

   Let    us   consider this time, the  linearization   at  the  non-grazing moment  $\omega =\pi.$  The   discontinuity    point  is    $c=(0, -\exp(0.0005\pi))$  and it is    of $(\alpha)-$ type, since $$\langle \nabla\Phi(c),f(c)\rangle=\langle(1,0)( -\exp(0.0005\pi), -0.001\exp(0.0005\pi))\rangle=-\exp(0.0005\pi)\neq 0.$$

%
%
%
%
%
%
%
%
%
%
%

By using formula (\ref{dert}), one can compute the gradient as $\displaystyle{\nabla\tau(c)=(\exp(-0.0005\pi),0)}.$ 
 
Then, utilizing  $\nabla\tau(c)$ and formula (\ref{derw}), one can determine   that   the  matrix   of   linearization  at  the moment $\pi$   is
 
\begin{eqnarray*}
B_2= \begin{bmatrix}
 \exp(-0.0005\pi)& 0\\
0.001& 0
   \end{bmatrix}.
\end{eqnarray*}

From the  monotonicity  of the jump  function, $-R_1y_2^2,$ it   follows that  the    the  yellow and blue subregions  
of $G$  are  invariant.
Consequently,  for   each solution near   to  $\Psi(t),$    the   sequences $B_i$ is of two   types 
$B_i=D_i^{(j)},$ $i\in\mathbb{Z}$ and $j=1,2,$  where $D_{2i-1}^{(1)}=O_{2},$  $D_{2i-1}^{(2)}=D_{1}^{(2)}=\begin{bmatrix} -1 & 0 \\  0.001-1.8\exp(0.00025\pi) &0\end{bmatrix},$  $D_{2i}^{(1)}=D_{2i}^{(2)}= \begin{bmatrix}
\exp(-0.0005\pi)& 0\\
0.001& 0
   \end{bmatrix}, i \in \mathbb Z.$   That  is, the   condition $(A4)$   is valid  and the linearization around the periodic solution (\ref{per}) on $\mathbb R$ is of two   subsystems:

\begin{equation}\label{orbexsub1}
\begin{aligned}
&u_1 '=u_2 , \\
&u_2 '=-u_1+0.001u_2,\\
&  \Delta u|_{t=\theta_{2i-1}}=  D_{2i-1}^{(1)}u,\\	
&  \Delta u|_{t=\theta_{2i}} =  D_{2i}^{(1)}u,\\	
\end{aligned}
\end{equation}

and  

\begin{equation}\label{orbexsub2}
\begin{aligned}
&u_1 '=u_2 , \\
&u_2 '=-u_1+0.001u_2,\\
&\Delta u|_{t=\theta_{2i-1}}=  D_{2i-1}^{(2)}u,\\	
&\Delta u|_{t=\theta_{2i}} =  D_{2i}^{(2)}u,\\	
\end{aligned}
\end{equation}
where   $\theta_{2i-1}=\frac{(2i-1)\pi}{2}$ and $\theta_{2i}=i\pi$
  
The sequences $\{D_{i}^{(j)}\},\ j=1,2,$ are $2-$ periodic.  It is appearant that system (\ref{orbexsub1})+(\ref{orbexsub2})  is a $(\omega,2)-$ periodic. 
Thus,   the  variational   system   for  the grazing  solution is   constructed.


\end{example}

\section{Orbital stability} \label{secorbitalstability}

\begin{figure}[ht] 
\centering
\includegraphics[width=5 cm]{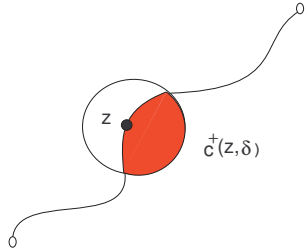}
\caption{The region $c^+(z,\delta).$}
\label{Cplus}
\end{figure}

In this section,  we    proceed   investigation of the  grazing   periodic solution $\Psi(t).$  Analysis  of  orbital stability  will be taken into account. Denote by $B(z ,\delta),$ an open ball with center at $z$ and the radius $\delta>0$ for a fixed point $z\in\Gamma\setminus\partial\Gamma.$   By condition (C3), the ball is divided by surface $\Gamma$ into two connected open regions. Denote $c^+(z,\delta),$ for the region, where solution  $x(t) = x(t,0,z)$ of \eqref{ode} enters as time increases. The region is depicted in Figure \ref{Cplus}. 

Set the  path of the periodic solution $\Psi(t)$   as  $$\eta:=\{x\in D: x=\Psi(t), \quad t \in \mathbb{R}\}.$$ 

%
%

%

Define  $ dist(A,a)=\inf_{\alpha \in A}\|\alpha-a\|,$  where  $A$  is a set, and $a$ is a point. 

\begin{definition}\label{orbitalstability}The   periodic solution $\Psi(t):\mathbb{R}\rightarrow D$ of (\ref{eq:graz1}) is said to be orbitally stable if for every $\epsilon>0,$ there corresponds $\delta = \delta(\epsilon)>0$ such that $dist(x(t,0,x_0),\eta)<\epsilon,$ for all $t\geq 0,$ provided $dist(x_0,\eta)<\delta $  and $x_0\notin \cup_i c^+(\Psi(\theta_i),\delta),\ for \ i=1,\ldots,m,$ where $m$ is the number of points $\Psi(\theta_i)\in\Gamma\setminus\partial \Gamma.$
\end{definition}

The  point     $x_0$  is not    considered in regions   $c^+(\Psi(\theta_i),\delta),i=1,\ldots,m,$  since solutions  which  start there   move   continuously  on  a finite interval,  while  $\Psi(t)$ experiences a non-zero jump at $t=\theta_i$ and this violates the continuity in initial value, in   general. In the same time, we take into   account any region adjoint to   points   of $\partial\Gamma,$ since the jump of $\Psi(t)$ is zero there  and, consequently,  the continuous   dependence in initial value is valid for all near  points.
 
\begin{definition}\label{asymptoticphase} The solution $\Psi(t):\mathbb{R_+}\rightarrow D$ of (\ref{eq:graz1}) is said to have asymptotic phase property if a $\delta>0$ exists such that to each $x_0$ satisfying $dist(x_0 ,\eta )<\delta$ and $x_0\notin \cup_i c^+(\Psi(\theta_i),\delta),\ for \ i=1,\ldots,m,$ there corresponds an asymptotic phase $\alpha(x_0) \in \mathbb R$ with property:  for all $\epsilon>0,$ there exists $T(\epsilon)>0,$  such that $x(t+\alpha(x_0),0,x_0)$ is in $\epsilon$-neighborhood of $\Psi(t)$ in $B-$topology for $t\in[T(\epsilon),\infty).$
\end{definition}


Let us consider the following system, which will be needed in the following lemmas and theorem 

\begin{equation}\label{orb}
\begin{aligned}
&x'=A(t)x, \\
&\Delta x|_{t=\zeta_i}= B_{i}u, 
\end{aligned}
\end{equation}
where  $A(t)$  and $B_i$   are $n\times n$  function-matrices,   $A(t+\omega)=A(t),$ for all $t\in\mathbb{R}$ and there exists an integer $p$ such that $\zeta_{i+p}=\zeta_i+\omega$ and $B_{i+p}=B_i,$ for all $i\in\mathbb{Z}.$

\begin{lemma}\label{Lemmafund}
Assume that system (\ref{orb}) has  a simple unit  characteristic multiplier and the remaining $n-1$ ones are in modulus less than unity. Then, the system (\ref{orb}) has a real fundamental matrix $X(t),$ of the form
\begin{equation}\label{fundamental1}
X(t)=P(t)\left( \begin{array}{cc}
1 & 0\\
0 & \exp{(Bt)} \end{array} \right), 
\end{equation}
where $P\in PC^1(\mathbb{R},\theta)$ is a regular, $\omega$-periodic matrix, and $B$ is an $(n-1)\times (n-1)$ matrix with all eigenvalues have negative real parts. 
\end{lemma}

\noindent \textbf{Proof.}
Denote the matrix $X(t), X(0)=I,$ as fundamental matrix of system (\ref{orb}). There   exists  a matrix $B_1$ such that   the   substitution  $x=P(t)z,$ where $P(t)=X(t)\exp(-B_1t),$   transforms (\ref{orb}) to the following system with constant coefficient \cite{ref1}, 

\begin{equation}\label{p-eq1}
\begin{aligned}
&z'=\Lambda z.
\end{aligned}
\end{equation}

The matrix  $ \exp(\Lambda \omega)$   has  a simple   unit eigenvalue and remaining $(n-1)$ ones   are in modulus less than unity.  Hence, there exists real nonsingular matrix $M,$  which  satisfies  
$$M^{-1}\exp(\Lambda \omega)M = \begin{bmatrix}
       1 & 0         \\[0.3em]
       0           & C_1 \\[0.3em]
      
     \end{bmatrix}.$$     
The remaining part of the proof is same as proof of Lemma 5.1.1 in \cite{ref1f}. 
     $\square$

Throughout  this section,    we will assume that $(A4)$ is valid.   That   is, the variational system (\ref{eq:graz1lin1}) consists of $m$ periodic subsystems.  For each of these systems, we find the matrix of monodromy, $U_j(\omega)$ and denote  corresponding  Floquet multipliers  by $\rho_i^{(j)},$ $i=1,\ldots,n,$  $j=1,\ldots,m.$  In the next part of the paper,  the following assumption is  needed.
\begin{itemize}
\item[(A5)]  $\rho_1^{(j)}=1$ and $|\rho_i^{(j)}|<1,$ $i=2,\ldots,n$ for each  $j=1,\ldots,m.$ 
\end{itemize}

\begin{lemma}\label{lemma2}
Assume that the assumptions (A4) and (A5) are valid. Then, for each  $j=1,\ldots,m,$  the system (\ref{eq:graz1lin1}) admits a fundamental matrix of the form
\begin{equation}\label{Fund1}
U_j(t)=P_j(t)[1,\exp(H_j \omega)], \quad t\in\mathbb{R},
\end{equation}
\end{lemma}
where $P_j\in PC^1(\mathbb{R},\zeta)$ is a regular, $\omega$-periodic matrix and $H_j$ is an $(n-1)\times (n-1)-$ matrix with all eigenvalues have negative real parts. 

The proof of Lemma \ref{lemma2}, can be done similar to that of Lemma \ref{Lemmafund}.


\begin{theorem}\label{thmorbitalstability} Assume that conditions $(C1)-(C7),$ $(C10),$ and the assumptions $(A1)-(A5)$ hold. Then $\omega$- periodic solution $\Psi(t)$ of $(\ref{eq:graz1})$ is orbitally asymptotically stable and has the asymptotic phase property.
\end{theorem}
\noindent \textbf{Proof.}  Since of the group property, we may assume $\Psi(0)$ is not a discontinuity point. Then, one can displace  the origin  to the point $\Psi(0),$ and the coordinate system can be rotated in such a way that the tangent vector $\Psi'_0=\Psi'(0)$ points in the direction of the positive $x_1$ axis i.e. the coordinates of this vector are $\Psi'_0=(\Psi'_{01},0,\ldots,0),$ $\Psi'_{01}>0.$

Let    $\theta_i, i\in\mathbb{Z},$    be  the  discontinuity moments of  $\Psi(t).$  Denote the path of   the solution by $\eta=\{x\in X: x=\Psi(t),t\in\mathbb{R}\}.$   There    exists a natural number  $p,$  such  that   $\theta_{i+p}=\theta_i +\omega$  for all  $i.$ Because of conditions $(C1)-(C7)$   and $K-$differentiability  of $\Psi(t)$   there exists continuous dependence on initial data and consequently there exists a neighborhood of $\eta$ such that any solutions which starts in the set  will have moments of discontinuity which constitute a $B-$ sequence with difference between neighbors approximately equal to the distance between corresponding neighbor moments of discontinuity of the periodic solution $\Psi(t).$ Consequently we can determine variational system for $\Psi(t),$ with points of discontinuity $\theta_i, \ i\in \mathbb{Z}.$

On the basis of discussion in Section  $2.1,$ one can define in the neighborhood of $\eta$ a $B-$ equivalent system of type (\ref{eq:graz2}).    The   variational    system of it  takes  the   form

\begin{equation}\label{orbital1}
\begin{aligned}
& z^\prime=A(t)z+r(t,z), \\
&\displaystyle{\Delta z|_{t=\theta_i}}= D_i^{(j)}z+q_i(z),\quad j=1,2,\ldots,m,
\end{aligned}
\end{equation}
where $r(t,z)=[f(\Psi(t)+z)-f(\Psi(t))]- A(t)z$ and $q_i(z)=W_i(\Psi(\theta_i)+z)-W_i(\Psi(\theta_i))-D_i^{(j)}z,$ are continuous functions,   and matrices  $D_i^{(j)}$  satisfy  condition $(A4).$  The functions are continuously differentiable with respect to $z.$  One can   verify that $r(t,0)\equiv q_i(0)\equiv 0$  and $r(t+\omega,z)=r(t,z)$ for $t\in\mathbb{R}.$   Moreover,   the derivatives  satisfy  $r'(t,0)\equiv q_{iz}'(0)\equiv 0$ and the functions  $r(t,z)\rightarrow 0,$ $q_i(z)\rightarrow 0,$ $r'_{z}(t,z)\rightarrow 0$ and $q_{iz}'(z)\rightarrow 0,$ as $z\rightarrow 0$ uniformly in $t\in[0,\infty), \ i\geq 0.$  Each system \eqref{orbital1} for  $j=1,2,\ldots,m,$ corresponds to a region adjoint to initial value, $x_0$ such that these regions cover a neighborhood of $x_0.$ 
%
%
%

Fix  a number $j$  and   denote  $Y_j(t)$  the fundamental matrix of  adjoint to  \eqref{orbital1}    linear homogeneous system
\begin{equation}\label{orbital1a}
\begin{aligned}
& y^\prime=A(t)y,\\
&\displaystyle{\Delta y|_{t=\theta_i}}= D_i^{(j)}y,
\end{aligned}
\end{equation}
of the form (\ref{Fund1}). One can verify that 
\begin{equation}\label{refleb}
Y_j(t)Y_j^{-1}(s)=P_j(t)\left( \begin{array}{cc}
1& 0\\
0 & \exp(H_j(t-s)) \end{array} \right)P_j^{-1}(s),
\end{equation}
for $-\infty<t,s<\infty .$ 

We can write 

\begin{eqnarray*}
\left( \begin{array}{cc}
1 & 0\\
0 & \exp(H_j(t-s)) \end{array} \right)=\left( \begin{array}{cc}
0 & 0\\
0 & \exp(H_j(t-s)) \end{array} \right)+\left( \begin{array}{cc}
1 & 0\\
0 & O_{n-1} \end{array} \right),
\end{eqnarray*}
where $O_{n-1}$ is the $(n-1)\times(n-1)$ zero matrix. Then it can be driven

$$Y_j(t)Y_j^{-1}(s)=G_1^{(j)}(t,s)+G_2^{(j)}(t,s)=G^{(j)}(t,s),$$

where

\begin{eqnarray*}
&G_1^{(j)}(t,s)=P_j(t)\left( \begin{array}{cc}
0& 0\\
0 & \exp(H_j(t-s)) \end{array} \right)P_j^{-1}(s),\\
&G_2^{(j)}(t,s)=P_j(t)\left( \begin{array}{cc}
1& 0\\
0 & O_{n-1} \end{array} \right)P_j^{-1}(s).
\end{eqnarray*}

Denote the eigenvalues of the matrix $H_j$ by $\lambda_2^{(j)}, \ldots, \lambda_n^{(j)}.$ By means of the Lemma \ref{Lemmafund} and \ref{lemma2},   there   exits  a number $\alpha>0,$ such that $Re(\lambda_k^{(j)})<-\alpha, $ $k=2,3,\ldots,n,$ where $Re(z)$ means the real part of the number, $z.$ Taking into account that the matrices $P_j$ and $P_j^{-1}$ are regular and periodic,  the following estimates can be calculated 
\begin{equation}\label{gest}
|G_1^{(j)}(t,s)|\leq K^{(j)}\exp(-\alpha(t-s)),
\end{equation}

\begin{equation}\label{gest2}
|G_2^{(j)}(t,s)|\leq K^{(j)},
\end{equation}

where $K^{(j)}$ is a positive real constant.

Denote the first column of the fundamental matrix $Y$ by $\chi^1.$ By the equation (\ref{Fund1}), $\chi^1$ is equal to the first column of $P_j,$ this means that it is a $\omega$-periodic solution of (\ref{eq:graz1lin1}). 

By assumptions of the theorem the variational system (\ref{orbital1}) satisfies the conditions of Lemma \ref{lemma2}, and  one can verify   that  the following estimate is true \cite{ref1f}

\begin{equation}\label{estimatefund}
|Y_j(t)|\leq K_1^{(j)}\exp(-\alpha t)\ \ for \ \ t\geq 0,
\end{equation}
where $K_1^{(j)}$ is a positive constant. Let us   setup  the following integral equation 
\begin{eqnarray}\label{estimatefund1}
&&z^{(j)}(t,a)=Y_j(t)a+\int\limits^{t}_{0}{G_1^{(j)}(t,s)r(s,z(s))ds}-\int\limits^{\infty}_{t}{G_2^{(j)}(t,s)r(s,z(s))ds}\nonumber\\
&&+\sum\limits_{0<\theta_k<t}{G_1^{(j)}(t,\theta_k+)q_k(z(\theta_k))}-\sum\limits_{t<\theta_k<\infty}{G_2^{(j)}(t,\theta_k+)q_k(z(\theta_k))},
\end{eqnarray}
where $a=[0,a_2,\ldots,a_n],$ $a_i\in\mathbb{R},$ $i=2,3,\ldots,n,$   are orthogonal to $\Psi^\prime(0),$ i.e. with the zero first coordinate.

Let $z_0^{(j)}(t,a)\equiv 0,$ and consider the following successive approximations

\begin{equation}\label{estimatefund2}
z_k^{(j)}(t,a)=Y_j(t)a+\int \limits^{\infty}_{0}{G^{(j)}(t,s)r(s,z_{k-1}(s))ds}+\sum\limits^{\infty}_{k=1}{G^{(j)}(t,\theta_k+)q_k(z_{k-1}(\theta_k))},
\end{equation}
for  $k =1,2,\ldots.$
By using the approximation (\ref{estimatefund2}) and estimation (\ref{estimatefund}),  one can verify that 
\begin{equation}\label{z1}
|z_1^{(j)}(t,a)|\leq K_1^{(j)}|a|\exp(-\alpha t/2).
\end{equation}
We will show that the bounded solution of  (\ref{estimatefund1}) exists  and satisfies \eqref{orbital1}.  For arbitrary positive small number  $L $, there exists a   number $ \delta =\delta(L)$ such that for $|z_1|<\delta,$  $|z_2|<\delta$  
\begin{equation}\label{lipsh1}
|r(t,z_1)-r(t,z_2)|\leq L|z_1-z_2|
\end{equation}
and
\begin{equation}\label{lipsh2}
|q_i(z_1)-q_i(z_2)|\leq L|z_1-z_2|,
\end{equation}
uniformly in $t\in[0,\infty).$

Denote by $\displaystyle{L_1=4K^{(j)}\Big(\frac{2}{\alpha}-\frac{1}{1-\exp(-\alpha\underline{\theta}/2)}\Big)}.$

Next, by using mathematical induction, we are going to show that $z_s^{(j)}(t,a), s=1,2,\ldots,$ are defined for $t\in[0,\infty)$ and satisfy 

\begin{equation}\label{approx}
|z_{s+1}^{(j)}(t,a)-z_{s}^{(j)}(t,a)|\leq K_1^{(j)}|a|\exp(-\alpha t/2)/2^s, \, s=0,1,2,\ldots,
\end{equation} 
if $L < L_1.$
Utilizing Lemma \ref{lemma2} and  inequalities (\ref{estimatefund}), (\ref{lipsh1}), (\ref{lipsh2})  and $\theta_{i+1}-\theta_i\geq \underline{\theta}, i \in \mathbb Z,$ one can verify that
\begin{equation}\label{realest}
|z_{k+1}^{(j)}(t,a)-z_{k}^{(j)}(t,a)|\leq K_1^{(j)}|a|L_1\exp(-\alpha t/2)/(2^{k}\alpha).
\end{equation}
As a consequence of \eqref{approx},  the sequence $z_{k+1}^{(j)}(t,a)$ converges uniformly on $t\in[0,\infty), \ |a|<\delta/2K_1^{(j)},$ and 
$$|z_{s}^{(j)}(t,a)|\leq 2K_1^{(j)}|a|\exp(-\alpha t/2), s = 1,2,\ldots.$$
Therefore, the  limit  function $z^{(j)}(t,a)$ exists on the same domain, it is piecewise continuous, satisfies  (\ref{estimatefund1}) and the following estimate
\begin{equation}\label{z}
|z^{(j)}(t,a)|\leq 2K_1^{(j)}|a|\exp(-\alpha t/2).
\end{equation}

Denote by $z(t)=z^{(j)}(t,a),$ for $j=1,2,\ldots,m.$
Next, we will verify that $z^{(j)}(t,a)$  satisfies (\ref{orbital1}). For it,  differentiate \eqref{estimatefund1}
\begin{eqnarray*}
&&z'(t)=Y^{\prime} _j(t)a+G_1^{(j)}(t,t)r(t,z(t))+G_2^{j}(t,t)r(t,z(t))+\int\limits^{t}_{0}{G_{1t}^{(j)}(t,s)r(s,z(s))ds}\\
&&-\int\limits^{\infty}_{t}{G_{2t}^{(j)}(t,s)r(s,z(s))ds}+\sum\limits_{0<\theta_k<t}{G_{1t}^{(j)}(t,\theta_k+)q_k(z(\theta_k))}-\sum\limits_{t<\theta_k<\infty}{G_{2t}^{(j)}(t,\theta_k+)q_k(z(\theta_k))}\\
&&=A(t)Y_j(t)a+G^{(j)}(t,t)r(t,z(t))+\int\limits^{\infty}_{0}{A(t)G^{(j)}(t,s)r(s,z(s))ds}+\sum\limits_{0<\theta_i<t}{A(t)G^{(j)}(t,\theta_k+)q_k(z(\theta_k))}\\
&&=A(t)z(t)+r(t,z(t)).
\end{eqnarray*}
 Fix  $\theta_k,$ $k\in\mathbb{Z},$ then

\begin{eqnarray*}
z(\theta_k+)-z(\theta_k)&&=Y_j(\theta_k+)a+\int\limits^{\theta_k}_{0}{G_1^{(j)}(\theta_k+,s)r(s,z(s))ds}-\int\limits^{\infty}_{\theta_k}{G_2^{(j)}(\theta_k+,s)r(s,z(s))ds}\\
&&+\sum\limits_{0\leq\theta_i<\theta_k}{G_1^{(j)}(\theta_k+,\theta_i+)q_i(z(\theta_i+))}-\sum\limits_{\theta_k<\theta_i<\infty}{G_2^{(j)}(\theta_k+,\theta_i+)q_i(z(\theta_i+))}\\
&&-Y_j(\theta_k)a-\int\limits^{\theta_k}_{0}{G_1^{(j)}(\theta_k,s)r(s,z(s))ds}+\int\limits^{\infty}_{\theta_k}{G_2^{(j)}(\theta_k,s)r(s,z(s))ds}\\
&&-\sum\limits_{0\leq\theta_i<\theta_k}{G_1^{(j)}(\theta_k,\theta_i+)q_i(z(\theta_i))}+\sum\limits_{\theta_k\leq\theta_i<\infty}{G_2^{(j)}(\theta_k,\theta_i+)q_i(z(\theta_i+))}\\
&&=D_{k}^{(j)}z(\theta_k)+q_k(z(\theta_k)).
\end{eqnarray*}
The above discussion proves  that $z^{(j)}(t,a), j = 1,,2,\ldots,m,$ are  bounded solutions of system (\ref{orbital1}).

We will determine the initial values of    bounded   solutions in terms of  $(n-1)$  parameters $a_2^{(j)},\ldots,a_n^{(j)}, j = 1,,2,\ldots,m.$   Denote  $a^{(j)}  = [0,a_2^j,a_3^j,\ldots,a_n^j].$ By using (\ref{estimatefund1}), we obtain

\begin{eqnarray*}
z^{(j)}(0,a^{(j)})&&=Y_j(0)a^{(j)}-\int\limits^{\infty}_{0}{G_{2}^{(j)}(0,s)r(s,z(s))ds}-\sum\limits_{0<\theta_k<\infty}{G_2^{(j)}(0,\theta_k+)q_k(z(\theta_k))}\\
&&=P_j(0)a^{(j)}-P_j(0)\left( \begin{array}{cc}
1& 0\\
0 & O_{n-1} \end{array} \right )\int\limits^{\infty}_{0}{P_j^{-1}(s)r(s,z(s))ds}-\sum\limits_{0<\theta_k<\infty}{P_j^{-1}(s)q_k(z(\theta_k))}.
\end{eqnarray*}
In the way utilized in \cite{ref1f}, one can   show  that  the coordinates of the initial value $(x_1,\ldots,x_n)\in D$   of the   solution  $z^{(j)}$  satisfy the equation 
\begin{equation}\label{surface}
x_1+\sum\limits_{i=2}^{n}{c_i^{j}x_i-h_j(x_2,\ldots,x_n)=0},
\end{equation}
where $h_j\in C^1, j = 1,,2,\ldots,m.$

One can see that equation (\ref{surface}) determines $(n-1)$ dimensional hypersurfaces $S^{j}\subset D,$ $j=1,2,\ldots,m,$ in a neighborhood of the origin such that each solution which starts at the surface satisfies inequality (\ref{z}). From the analytical representation, it follows  that the equation of the tangent space of $S^j$ at the origin  is described by the equation $x_1+\sum\limits_{i=2}^{n}{c_i^{j}x_i}$ and the first coordinate of the gradient of the left hand side in (\ref{surface}) is unity.  Moreover,  the path $\eta$ intersects $S^j$ transversely. This   and condition   $(A4)$   imply    that the path of every solution  $\phi(t)$   near  $\Psi(t)$ intersects    one of the manifolds $S^j, j = 1,2,\ldots,m,$ at some $\bar t\in[0,2\omega].$

Because of the continuous dependence on initial values, a $\delta(\epsilon)>0$ exists for a given $\epsilon>0,$ such that if $dist(x^0,\eta_{\delta})<\delta(\epsilon),$ then the solution $\phi(t,x^0)$ is defined on $[0,2\omega],$ and $dist(\phi(t,x^0),\eta)<\epsilon\leq \epsilon_1$  for $t\in[0,2T].$ Therefore, the path of $\phi(t,x^0)$ intersects $S^j$ for some $j=1,2,\ldots,m$ and $t_1\in[0,2\omega].$ The solution $\phi(t,\phi(t_1,x^0))=\phi(t+t_1,x^0)$ has its initial value in $S^j,$ consequently, satisfies (\ref{z}).  In the light of the $B-$ equivalence, the corresponding solution $x(t), x(0)=\phi(0)-\Psi(0),$ of \eqref{orbital1} satisfies the property that for all $\epsilon>0,$ there exists $T(\epsilon)$ such that $x(t)$ is in an $\epsilon-$ neighborhood of $\Psi(t)$ for $t\in [T(\epsilon),\infty).$ That is,  the solution $\Psi(t)$  is orbitally asymptotically stable and there exists an asymptotical phase. $\square$

Definitions of the orbital stability and an asymptotic phase as well as theorem of orbital stability for non-grazing periodic solutions are  also  presented in \cite{simo-bain}. In our paper, we suggest the orbital stability theorem for grazing periodic solutions, its proof and formulate the definitions for the stability. They are different in many aspects from those provided in \cite{simo-bain}. It   is valuable that   they  also  valid, if the solution is non-grazing.

To shed light on our theoretical results, we will present the following examples.

\begin{example}\label{ex1}  We continue with the system presented in Examples \ref{exdds}  and \ref{linear}.
In Example \ref{exdds}, we verified that system (\ref{orbex22})defines a $K-$ smooth discontinuous flow    in the   plane and   the variational system  (\ref{orbexsub1})+(\ref{orbexsub2}) around the  grazing  periodic   solution, $\Psi(t)$ is approved.

 Using systems (\ref{orbexsub1})  and (\ref{orbexsub2}), one can evaluate the Floquet multipliers  as  $\rho_1^{(1)}=1,$  $\rho_2^{(1)}=0.8551,$ $\rho_1^{(2)}=1$ and $\rho_2^{(2)}=0.$ This verifies condition $(A5).$ 

The conditions $(C1)-(C7)$ and $(C10)$ are validated and the assumptions (A4) and (A5) verified. By using Theorem \ref{thmorbitalstability}, we can assert that the solution, $\Psi(t)$  is orbitally asymptotically stable. The stability is illustrated in   Fig. \ref{disc}. The red one is for a trajectory of  the discontinuous periodic solution (\ref{per}) of (\ref{orbex2}) and the blue one is for the  near  solution of (\ref{orbex2}) with initial value $y_0=(0.8,1.2).$ It can be observed from Fig. \ref{disc} that the blue trajectory approaches the red one as time increases. 

 \begin{figure}[htbp] 
\centering
\includegraphics[width=9 cm]{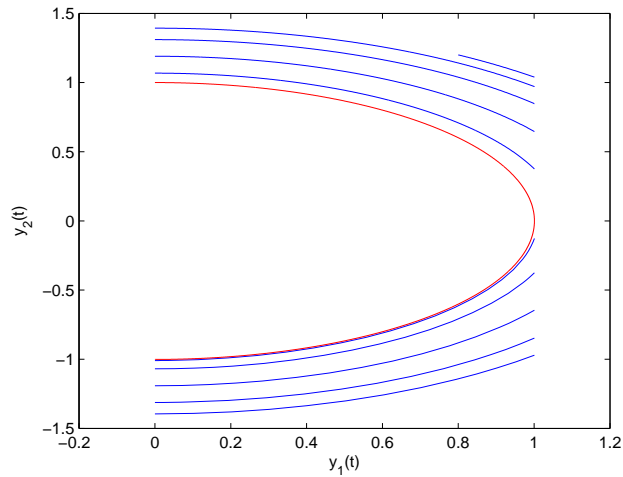}
\caption{The red   discontinuous  cycle of (\ref{orbex2}) axially  grazes $\Gamma$  at $(0.00025\pi),0)$  and     $(0, -\exp(-0.0005\pi))$  is an  $(\alpha)$-type point.  The blue arcs are of the  trajectory with initial value $(0.8, 1.2).$  It can be observed  that   it approaches the   grazing one as time increases.}
\label{disc}
\end{figure}

\end{example}

\begin{example} (A   periodic   solution with a non-axial grazing).  We will take into account the following autonomous system with variable moments of impulses 

\begin{equation}\label{exnonaxial}
\begin{aligned}
&x_1'=x_2,\\
&x_2'= -x_1,\\
&\displaystyle{\Delta x_1|_{x\in \Gamma }= \frac{1}{\sqrt{2}}-x_1+K(x_2-x_1)^2,}\\
&\displaystyle{\Delta x_2|_{x\in \Gamma }=  \frac{1}{\sqrt{2}}-x_2+K(x_2-x_1)^2,}
\end{aligned}
\end{equation}
where $\Gamma=\{(x_1,x_2)| x_1+x_2=\sqrt{2}\},$ $\tilde \Gamma=\{(x_1,x_2)| x_1=x_2\}$ and $K=0.11.$ It is easy to verify that the point $\displaystyle{x^*=(\frac{1}{2},\frac{1}{2})}$ is a grazing point because $\langle\nabla\tau(x^*),f(x^*)\rangle=\langle(1,1),(\frac{1}{2},-\frac{1}{2})\rangle=0$  and the grazing  is non-axial. We   assume that  the domain is  the plane. 

The solution $\Psi(t)=(\sin(t),\cos(t)),$ $t\in\mathbb R$ is a grazing   one,   since the  point $x^*=\Psi(\frac{\pi}{4})$ is  from   its orbit. The    cycle and the line of discontinuity  are depicted in Figure \ref{nonaxialcircle}.

\begin{figure}[ht] 
\centering
\includegraphics[width=9 cm]{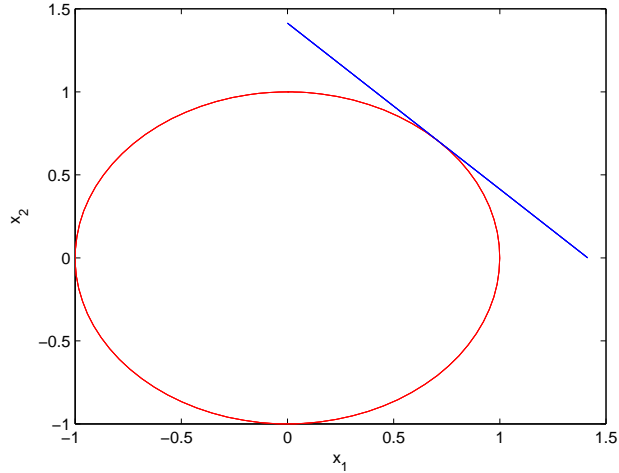}
\caption{The red curve is  the  orbit  of  $\Psi(t)$   which grazes  non-axially the   line of discontinuity.}
\label{nonaxialcircle}
\end{figure}

Let us consider the linearization at the grazing point $x^*$   next. We will consider the near solution $x(t)=x(t,0,x^*+\Delta x).$   Denote $t=\xi,$    the moment when  the solution meets the surface of discontinuity $\Gamma$ at the point $\bar x = x(\xi)=x(\xi,0,x^*+\Delta x).$  Taking into account formulae \eqref{derW_1}, \eqref{derW_2} with \eqref{exnonaxial},   one  can obtain the following matrix

\begin{eqnarray}\label{derW_222}
&\displaystyle{\frac{\partial W_i(x(\xi,0,x^*+\Delta x))}{\partial x_1^0}}=\begin{bmatrix} \bar x_2\\
 -\bar x_1 \end{bmatrix}\frac{1}{\bar x_1-\bar x_2}\displaystyle{+\begin{bmatrix} -2K(\bar x_2-\bar x_1) &  -2K(\bar x_2-\bar x_1) \\ -2K(\bar x_2-\bar x_1) &  -2K(\bar x_2-\bar x_1)\end{bmatrix}} \nonumber \\
&\times\Bigg(e_1+\begin{bmatrix} \bar x_2\\
 -\bar x_1 \end{bmatrix}\frac{1}{\bar x_1-\bar x_2}\Bigg)+\displaystyle{\begin{bmatrix} -\frac{1}{\sqrt{2}}+K(\bar x_2-\bar x_1)^2\\
 \frac{1}{\sqrt{2}}-K(\bar x_2-\bar x_1)^2 \end{bmatrix}\frac{1}{\bar x_1-\bar x_2}}. 
\end{eqnarray}

Calculating the right hand side of the expression \eqref{derW_222}, we obtain that

\begin{eqnarray}\label{derW_333}
&\displaystyle{\frac{\partial W_i(x(\xi,0,x^*+\Delta x))}{\partial x_1^0}}=\begin{bmatrix} \displaystyle{\frac{-\sqrt{2}+0.22}{ 2\sqrt{2}}} \\[6pt] \displaystyle{\frac{\sqrt{2}+0.22}{ 2\sqrt{2}}} \end{bmatrix}.
\end{eqnarray}

Using similar method with that of the first one, the second derivative can be computed as 

\begin{eqnarray}\label{derW_444}
&\displaystyle{\frac{\partial W_i(x(\xi,0,x^*+\Delta x))}{\partial x_2^0}}=\begin{bmatrix} \displaystyle{\frac{\sqrt{2}+0.22}{ 2\sqrt{2}}}  \\[6pt] \displaystyle{\frac{-\sqrt{2}+0.22}{ 2\sqrt{2}}}  \end{bmatrix}.
\end{eqnarray}

Combining \eqref{derW_333} and \eqref{derW_444}, we can obtain the following matrix for the linearization at the grazing point $x^*,$ 
\begin{eqnarray}\label{derW_555}
&\displaystyle{  W_{ix}(x^*)}=\displaystyle{ \begin{bmatrix} \displaystyle{\frac{-\sqrt{2}+0.22}{ 2\sqrt{2}}}  &\displaystyle{\frac{\sqrt{2}+0.22}{ 2\sqrt{2}}}  \\ \displaystyle{\frac{\sqrt{2}+0.22}{ 2\sqrt{2}}} &\displaystyle{\frac{-\sqrt{2}+0.22}{ 2\sqrt{2}}} \end{bmatrix}.}
\end{eqnarray}
It is appearant that the matrix $W_{ix}(x^*)$ is continuous with respect to its arguments,  since  it  is   constant  if the point $x^*+\Delta x$   is not   from the  orbit of the grazing  solution. 
Since   of the  limit  procedure,   it is the  same constant  for  all  points of the grazing  solution.  Thus,  the  Jacobian is constant matrix in a neighborhood of the  grazing  point  and  condition $(A2)$ is valid.

Now, let us check the validity of the condition $(A3).$ Consider a near solution $x(t)=x(t,0,\bar x),$ to the grazing cycle $\Psi(t),$ where $\bar x=(0,\bar x_2),\, \bar x_2>\Psi_2(0)=1.$ So, the near solution $x(t)$ satisfies the condition $(N1).$ For the grazing periodic solution, it is true that $\theta_{i+1}-\theta_i=2\pi=\omega.$  The grazing solution $\Psi(t) = x(t,0,(0,1)),$ touches the line of discontinuity $\Gamma$ at $t=\frac{\omega}{8}.$ The  first coordinate of the near solution   is  $x_1(t)=\bar x_2\sin(t),$  and  $\displaystyle{x_1(\frac{\omega}{8})=\bar x_2\sin(\frac{\omega}{8})=\frac{\bar x_2}{\sqrt{2}}>\Psi_1(\frac{\omega}{8})=\frac{1}{\sqrt{2}}}.$ Consequently,  the near solution $x(t)$ meets the line of discontinuity $\Gamma$ before the moment $\frac{\omega}{8}.$ This implies that $0<\tau(x)<\frac{\pi}{4}-\epsilon,$ for a small positive $\epsilon$ whenever  $x_1(t)$ is close to $\frac{1}{\sqrt{2}}.$ Thus, the condition $(A3)$ is valid and  Lemma \ref{conttau} proves  condition $(C).$


In the light of  the above discussion,  the  bivalued  matrix   of coefficients for the grazing point is easily obtained as

\begin{eqnarray}\label{B1}
 B_1=\begin{cases} O_{2}, \ & \mbox{if (N1) is valid}, \\ 
\begin{bmatrix} \displaystyle{\frac{-\sqrt{2}+0.22}{ 2\sqrt{2}}}  &\displaystyle{\frac{\sqrt{2}+0.22}{ 2\sqrt{2}}}  \\ \displaystyle{\frac{\sqrt{2}+0.22}{ 2\sqrt{2}}} &\displaystyle{\frac{-\sqrt{2}+0.22}{ 2\sqrt{2}}} \end{bmatrix},\ & \mbox{if (N2) is valid}. \end{cases}
\end{eqnarray}

It is appearant that the interior  of the   grazing orbit is invariant. Let us show that the external part of the unit circle is  positively   invariant. It is sufficient to demonstrate that $J_1(x_1)^2+J_2(x_2)^2>1$ for any $(x_1,x_2)\in \Gamma.$ Denote $x_1=z$ and $x_2=\sqrt{2}-z$ and consider the formula  $$F(z)= J_1(z)^2+J_2(\sqrt{2}-z)^2  = (\frac{1}{\sqrt{2}}+0.11(\sqrt{2}-2z)^2)^2+(\frac{1}{\sqrt{2}}+0.11(\sqrt{2}-2z)^2)^2,$$
where $F(\frac{1}{\sqrt{2}})=1.$ It is easy to calculate that $F'(\frac{1}{\sqrt{2}})=0$ and $F''(\frac{1}{\sqrt{2}})=0.88\sqrt{2}>0.$ Consequently, $$F(z)-F(\frac{1}{\sqrt{2}})=F(z)-1=\frac{1}{2}F''(\frac{1}{\sqrt{2}})(z-\frac{1}{\sqrt{2}})^2+o(\|z-\frac{1}{\sqrt{2}}\|^2)>0,$$ if $z$ is close to $\frac{1}{\sqrt{2}}.$ Thus, near the grazing point, the external region is invariant. From this discussion, since of the formula \eqref{B1}, we can conclude that  the condition  $(A4)$ is valid. Taking   into account it with the expression \eqref{B1}, the linearization system for \eqref{exnonaxial} around  the grazing solution  $\Psi(t)$ is obtained as
\begin{equation}\label{exnonaxiallin}
\begin{aligned}
&u_1'=u_2,\\
&u_2'= -u_1,\\
&\displaystyle{\Delta u(2\pi i)}=D_i^{(j)}u,
\end{aligned}
\end{equation}
where $D_i^{(1)}=O_2$ and $D_i^{(2)}=\begin{bmatrix} \displaystyle{\frac{-\sqrt{2}+0.22}{ 2\sqrt{2}}}  &\displaystyle{\frac{\sqrt{2}+0.22}{ 2\sqrt{2}}}  \\ \displaystyle{\frac{\sqrt{2}+0.22}{ 2\sqrt{2}}} &\displaystyle{\frac{-\sqrt{2}+0.22}{ 2\sqrt{2}}} \end{bmatrix}, \, i\in\mathbb{Z}.$ 

To finalize stability analysis, consider the first system in \eqref{exnonaxiallin}, with matrices $D_i^{(1)} = O_2.$ Its multipliers are $\rho_1^{(1)}=\rho_2^{(1)}=1$ and it constitutes the linearization for the orbits which are inside the circle. The system does not give a decision by orbital stability theorem, Theorem \ref{thmorbitalstability}. Nevertheless, from the simple analysis \cite{perko}  result, we know that the grazing orbit is stable with respect to inside orbits of the system.   The linearization  of  orbits which   are outside   of the circle   has multipliers $\rho_1^{(2)}=1$ and $\rho_2^{(2)}=-0.15.$ It means that the periodic solution is   orbitally  stable with respect to solutions outside of the circle. Summarizing the discussion, we can conclude that the periodic solution is stable. The   stability   result  is   observed through   simulations and it  is  seen  in  Fig. \ref{nonaxialcircle1}.

\begin{figure}[ht] 
\centering
\includegraphics[width=9 cm]{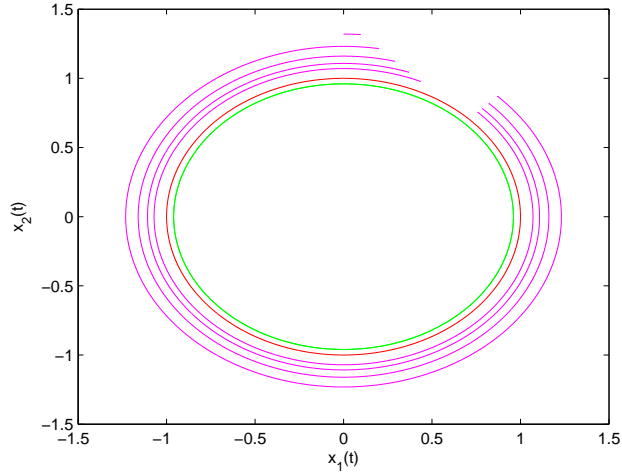}
\caption{The red  orbit   of  system \eqref{exnonaxial}   non-axially  grazes the   surface  $\Gamma.$   The magenta trajectory  with initial point $(0,1.32)$  approaches the cycle as time increases. The green cycle with initial point  $(0,0.96)$ demonstrates  the inside stability  of the   grazing   orbit.}
\label{nonaxialcircle1}
\end{figure}

\end{example}

\section{Small parameter analysis and  grazing bifurcation}
In this part, we will discuss  existence   and  bifurcation  of   cycles   for  perturbed systems,   if the generating   one  admits a grazing  periodic solution. In continuous dynamical  systems, a small parameter may cause a change in the number of periodic solutions  in critical cases.  In  the present   analysis,  we  will   demonstrate   that  the change  may  happen  in \textit{non-critical }  cases, since of the non-transversality.  That  is why,  one can  say  that   \textit{grazing  bifurcation}  is under  discussion.   
Let us  deal with the following  system

\begin{equation}\label{smallpar}
\begin{aligned}
& x'=f(x)+\mu g(x,\mu),\\
& \Delta x|_{x\in\Gamma(\mu)}=I(x)+\mu K(x,\mu),
\end{aligned}
\end{equation}
where $x\in  \mathbb R^n, t \in \mathbb R, \Gamma(\mu)=\{x| \ \Phi(x)+\mu\phi(x,\mu)=0\},$ $\mu\in(-\mu_0,\mu_0),$  and  $\mu_0$ is a sufficiently small positive number.    Functions  $f(x), I(x)$    and $\Phi(x)$   are  continuously differentiable up to second order,  $g(x,\mu),K(x,\mu) $  are  continuously differentiable in $x$ and $\mu.$     The   function  $\phi(x,\mu)$ is continuously differentiable in $x$ up to second order and  to  first order in $\mu. $ We assume that the generating system for \eqref{smallpar} is the system \eqref{eq:graz1} with all conditions assumed   for   the system,  earlier.  The   main assumption of this section  is that (\ref{eq:graz1}) admits a $\omega-$periodic solution, $\Psi(t).$    Let    $\Psi(0)=(\zeta_1^0,\zeta_2^0,\ldots,\zeta_n^0)$  be   the  initial value  of the   solution. 
 
 Our aim is to find conditions that verify the existence of periodic solutions of  \eqref{smallpar} with a period $\mathscr{T} $ such that for $\mu=0,$ the periodic solutions of \eqref{smallpar} are turned down to $\Psi(t).$ It is common for the autonomous systems that the period $\mathscr{T} $ does not coincide with $\omega.$ Thus, in the remaining part of the paper, we will consider the period $\mathscr{T} $ as an unknown  variable.

Since $\Psi(0)$ is not an equilibrium, there is a number $j=1,2,\ldots,n,$ such that $f_j(\zeta_1^0,\zeta_1^0,\ldots,\zeta_n^0)\neq 0.$ In other words, the vector field  is transversal to line $x_j=\zeta_j^0$  near the point.  Hence, to try points near to $\Psi(0)$ for the periodicity, it is sufficient to consider those with $j-$th coordinate is equal to $\zeta_j^0,$ 
 \cite{malkin}.  For the discontinuous dynamics, the choice of the fixed coordinate  can  be made  easier  if  the surface of discontinuity   is  provided with a constant coordinate.   We  will demonstrate this  in examples. 
Denote the initial values of the intended periodic solution by $\zeta_1,\zeta_2,\ldots,\zeta_n.$ Assume that one initial value $\zeta_j$ is known, i.e. $\zeta_j^0.$ Thus, the problem contains $n-$many unknowns, they can be presented as $\zeta_1,\zeta_2,\ldots,\zeta_{j-1},\zeta_{j+1},\ldots,\zeta_n,\mathscr{T}.$ Denote the solution of \eqref{smallpar} by $x_s(t,\zeta_1,\zeta_2,\ldots,\zeta_n,\mu)$ with initial conditions $x_s(0,\zeta_1,\zeta_2,\ldots,\zeta_n,\mu)=\zeta_s.$ To determine the unknowns, we will consider the Poincar$\acute{e}$  criterion, which can be written as
\begin{equation}\label{periodpar}
\begin{aligned}
\mathscr{S}_k(\mathscr{T},\zeta_1,\zeta_2,\ldots,\zeta_n,\mu) \equiv x_k(\mathscr{T},\zeta_1,\zeta_2,\ldots,\zeta_n,\mu)-\zeta_k=0, \, k=1,2,\ldots,n,
\end{aligned}
\end{equation}
where $\zeta_j=\zeta_j^0.$ The  equations \eqref{periodpar}    are   satisfied  with   $\mu=0, \mathscr{T}= \omega,  \zeta_i = \zeta_i^0, i = 1,2,\ldots,n,$  since  $\Psi(t)$  is  the periodic solution.  
  
%
  
The following condition for the determinant is also needed in the remaining part paper.

\begin{itemize}
\item[(A6)]   
\begin{equation}
\begin{vmatrix}
 \frac{\partial(\mathscr{S}_1(\omega,\zeta_1^0,\zeta_2^0,\ldots,\zeta_{j-1}^0,\zeta_{j+1}^0,\ldots,\zeta_n^0,0))}{\partial\mathscr{T}} & \ldots & \frac{\partial(\mathscr{S}_1(\omega,\zeta_1^0,\zeta_2^0,\ldots,\zeta_{j-1}^0,\zeta_{j+1}^0,\ldots,\zeta_n^0,0))}{\partial \zeta_n}\\ \frac{\partial(\mathscr{S}_2(\omega,\zeta_1^0,\zeta_2^0,\ldots,\zeta_{j-1}^0,\zeta_{j+1}^0,\ldots,\zeta_n^0,0))}{\partial\mathscr{T}} & \ldots & \frac{\partial(\mathscr{S}_2(\omega,\zeta_1^0,\zeta_2^0,\ldots,\zeta_{j-1}^0,\zeta_{j+1}^0,\ldots,\zeta_n^0,0))}{\partial \zeta_n} \\ \vdots & \ddots &\vdots\\ \frac{\partial(\mathscr{S}_n(\omega,\zeta_1^0,\zeta_2^0,\ldots,\zeta_{j-1}^0,\zeta_{j+1}^0,\ldots,\zeta_n^0,0))}{\partial\mathscr{T}} & \ldots & \frac{\partial(\mathscr{S}_n(\omega,\zeta_1^0,\zeta_2^0,\ldots,\zeta_{j-1}^0,\zeta_{j+1}^0,\ldots,\zeta_n^0,0))}{\partial \zeta_n} \end{vmatrix}\neq 0
\end{equation}
\end{itemize}

\begin{theorem}\label{theosmall1} Assume that  condition $(A6)$ is valid.  Then,  (\ref{smallpar}) admits a non-trivial periodic solution, which converges in the $B-$ topology to the non-trivial $\omega$-periodic solution of (\ref{smallpar}) as $\mu$ tends to zero. 
\end{theorem}

We will present the following examples to realize our theoretical results.

\begin{example} In this example, we will consider the perturbed system in case the generating system has a graziness. To show that, let us take into account the following perturbed system

\begin{equation}\label{smalldisc}
\begin{aligned}
&x_1'=x_2,\\
&x_2'=-0.001x_2-x_1,\\
&\Delta x_2|_{x\in \Gamma_1}= -(1+R_1x_2+\mu x_2)x_2,\\
&\Delta x_2|_{x\in \Gamma_2}= -(1+R_2+\mu (x_2-\exp(0.001\pi/2))x_2.
\end{aligned}
\end{equation}
It is easy to see that the system (\ref{smalldisc}) is of the form (\ref{smallpar}). For $\mu=0,$ the generating system became (\ref{orbex2}). For the perturbed system (\ref{smalldisc}), we will investigate existence of the  periodic solution around the grazing periodic solution of  (\ref{orbex2})  with  the  help   of Theorem \ref{theosmall1}.

There   are   two  sorts   of possible periodic  solutions of (\ref{smalldisc}) around the grazing one.  One of them has two impulse moments during the period since it crosses both lines of discontinuity, i.e. $x_1=0$ and $x_1=\exp(0.00025\pi).$ The other sort is the periodic solution  which  does not intersect the line $x_1=\exp(0.00025\pi) $ and intersects the line $x_1=0.$ We will show the existence of both type of periodic solutions if $|\mu|$ sufficiently small. 

Let us start with the second type, assume that the solution for the perturbed system exists and it starts   at   the point  $(0,x_{02}), x_{02}<1 $  and  does not intersect the line $x_1=\exp(0.00025\pi).$ Denote the initial values of the periodic solution by $\zeta_1$ and $\zeta_2.$ Since the periodic solution necessarily    intersects the line $x_1=0,$   one can      choose   $\zeta_1 \equiv \zeta_1^0=0.$    By specifying  the formula in \eqref{periodpar} for the system \eqref{smalldisc}, it is easy to obtain the following expressions
\begin{equation}\label{example1small}
\begin{aligned}
&\mathscr{S}_1(\mathscr{T},0,\zeta_2,\mu)= x_1(\mathscr{T},0,\zeta_2,\mu)=0,\\
&\mathscr{S}_2(\mathscr{T},0,\zeta_2,\mu)=x_2(\mathscr{T},0,\zeta_2,\mu)-\zeta_2=0.
\end{aligned}
\end{equation}

Next, taking the derivative of the expressions in \eqref{example1small}, we can obtain the following 
 
\begin{eqnarray}\label{smallexone1}
 \begin{vmatrix} \frac{\partial (\mathscr{S}_1(\mathscr{T},0,\zeta_2,\mu))}{\partial \mathscr{T} }&  \frac{\partial (\mathscr{S}_1(\mathscr{T},0,\zeta_2,\mu))}{\partial \zeta_2 }\\ \frac{\partial (\mathscr{S}_2(\mathscr{T},0,\zeta_2,\mu))}{\partial \mathscr{T} }&  \frac{\partial (\mathscr{S}_2(\mathscr{T},0,\zeta_2,\mu))}{\partial \zeta_2 } \end{vmatrix} =
 \begin{vmatrix}  \frac{\partial x_1(\omega,0,\zeta_2^0,0)}{\partial \mathscr{T}} & \frac{\partial x_1(\omega,0,\zeta_2^0,0)}{\partial \zeta_2}\\ \frac{\partial x_2(\omega,0,\zeta_2^0,0)}{\partial \mathscr{T}} & \frac{\partial x_2(\omega,0,\zeta_2^0,0)}{\partial \zeta_2}-1\end{vmatrix} .
\end{eqnarray}

The determinant \eqref{smallexone1} is calculated by means of  the monodromy matrix  of (\ref{orbex2}), with the impulse matrix  $D_1^{(1)}=O_2,$ i.e. 

\begin{equation}\label{monodromy11}
\begin{bmatrix} 1 &  -0.0317 \\ 1.0158 & -0.1014
\end{bmatrix}.
\end{equation}

Taking into account the system \eqref{smallexone} with \eqref{monodromy11} at $\zeta_2=\zeta_2^0$ and $\mathscr{T}=\omega$ for $\mu=0,$ one can derive that 

\begin{eqnarray}\label{smallexoneone}
 \begin{vmatrix} \frac{\partial \mathscr{S}_1(\omega,0,\zeta_2^0,0)}{\partial \mathscr{T} }&  \frac{\partial \mathscr{S}_1(\omega,0,\zeta_2^0,0)}{\partial \zeta_2 }\\ \frac{\partial \mathscr{S}_2(\omega,0,\zeta_2^0,0)}{\partial \mathscr{T} }&  \frac{\partial \mathscr{S}_2(\omega,0,\zeta_2^0,0)}{\partial \zeta_2 }-1 \end{vmatrix} =-0.0317\exp(0.00025\pi)\neq0.
\end{eqnarray}

This verifies condition $(A6).$  Thus, condition $(A6)$ is valid, then by utilizing Theorem \eqref{theosmall1},   we can assert that the system  (\ref{smallpar}) admits a non-trivial periodic solution, which converges in the $B-$ topology to the non-trivial $\omega$-periodic solution of (\ref{eq:graz1}) as $\mu$ tends to zero. 


Now, let us verify that system (\ref{smalldisc}) has a circle  which intersects the line $x_1=\exp(0.00025\pi)$ in the neighborhood of $(\exp(0.00025\pi),0).$ So, the periodic solution will attain two discontinuity moments in a period.  Denote the initial values of the periodic solution by $\zeta_1$ and $\zeta_2.$ To apply the condition $(A6),$ fix  one initial value $\zeta_1=\zeta_1^0=0$ of the intended periodic solution and in the light of the expressions \eqref{periodpar}

\begin{equation}\label{example1small2}
\begin{aligned}
&\mathscr{S}_1(\mathscr{T}, 0, \zeta_2,\mu)=x_1(\mathscr{T}, 0, \zeta_2,\mu)=0,\\
&\mathscr{S}_2(\mathscr{T}, 0, \zeta_2,\mu)=x_2(\mathscr{T}, 0,\zeta_2,\mu)-\zeta_2=0.
\end{aligned}
\end{equation}

Taking the derivative of the expressions \eqref{example1small2} with respect to variables $\mathscr{T}$ and $\zeta_2,$ one can obtain the following 

\begin{eqnarray}\label{smallexone}
 \begin{vmatrix} \frac{\partial (\mathscr{S}_2(\mathscr{T},0,\zeta_2,\mu))}{\partial \mathscr{T} }&  \frac{\partial (\mathscr{S}_1(\mathscr{T},0,\zeta_2,\mu))}{\partial \zeta_2 }\\ \frac{\partial (\mathscr{S}_2(\mathscr{T},0,\zeta_2,\mu))}{\partial \mathscr{T} }&  \frac{\partial (\mathscr{S}_2(\mathscr{T},0,\zeta_2,\mu))}{\partial \zeta_2 } \end{vmatrix} =
 \begin{vmatrix}  \frac{\partial x_1(\omega,0,\zeta_2^0,0)}{\partial \mathscr{T}} & \frac{\partial x_1(\omega,0,\zeta_2^0,0)}{\partial \zeta_2}\\ \frac{\partial x_2(\omega,0,\zeta_2^0,0)}{\partial \mathscr{T}} & \frac{\partial x_2(\omega,0,\zeta_2^0,0)}{\partial \zeta_2}-1\end{vmatrix} .
\end{eqnarray}

To determine the above determinant, the monodromy matrix  of (\ref{eq:graz1lin1}) with the jump matrix $D_i^{(2)}$ can be evaluated as

\begin{equation}\label{monodromy111}
 \begin{bmatrix} 1 & 0.01 \\ 0 & 0.704
\end{bmatrix} .
\end{equation}

For $\mu=0,$ with the values $\omega$ and $\zeta_2^0$ the determinant \eqref{smallexone} can be determined as 

\begin{eqnarray}\label{smallexoneoneoe}
 \begin{vmatrix} \frac{\partial (\mathscr{S}_1(\omega,0,\zeta_2^0,0))}{\partial \mathscr{T} }&  \frac{\partial (\mathscr{S}_1(\omega,0,\zeta_2^0,0))}{\partial \zeta_2 }\\ \frac{\partial (\mathscr{S}_2(\omega,0,\zeta_2^0,0))}{\partial \mathscr{T} }&  \frac{\partial (\mathscr{S}_2(\omega,0,\zeta_2^0,0))}{\partial \zeta_2 } \end{vmatrix}=
\begin{vmatrix}  0 & 0.01\\ -\exp(0.00025\pi) & -0.296 \end{vmatrix} =0.01\exp(0.00025\pi) \neq 0.
\end{eqnarray} 

This verifies condition $(A6).$ So, By Theorem \ref{theosmall1}, we can conclude that the perturbed system \eqref{smalldisc}   admits a non-trivial $\mathscr{T}(\mu)-$ periodic solution which converges in the $B-$ topology to the non-trivial $\omega$-periodic solution of (\ref{orbex2}) as $\mu$ tends to zero such that $\mathscr{T}(0)=\omega.$

%
%
%
%
%

In  Fig. \ref{figsmall}, some numerical results are provided to show the solutions of system (\ref{smalldisc}) with $\mu=0.05.$

\begin{figure}[ht] 
\centering
\includegraphics[width=9 cm]{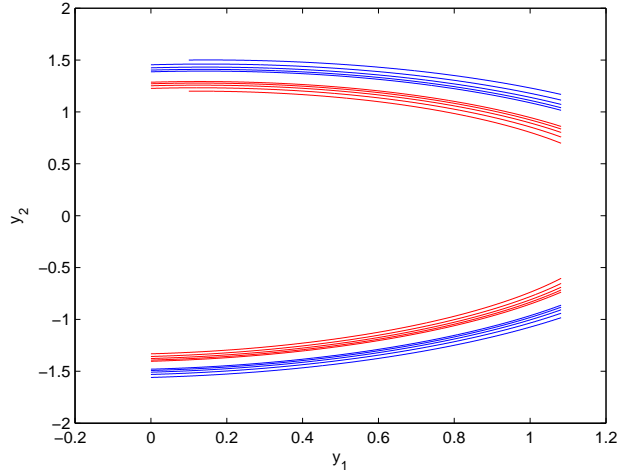}
\caption{The red arcs are the trajectory of the system (\ref{smalldisc}) with initial value $(0, 1.2)$ and the blue arcs are the orbit with initial value $(0, 1.5).$ Through simulation, we observe that the trajectories approach to the periodic solution of (\ref{smalldisc}) as time increases.}
\label{figsmall}
\end{figure}

The periodic solutions for $\mu\neq 0$ are not grazing.  For $\mu=0,$ we have one periodic solution which is orbitally stable, and for $\mu<0,$ there exist two  periodic solutions. One of them has  one discontinuity moment in each period, in other words, the cycle does not intersect the surface of discontinuity around grazing point and it is orbitally stable and the other one has two discontinuity moments in each period. This means, the number of periodic solutions increases by variation of $\mu,$   around $\mu=0.$ So, we will call that bifurcation of periodic solution from a grazing cycle.

\end{example}

\begin{example}\label{smallex1dof2}  Let us  consider  the following system with variable moments of impulses and  a small parameter
\begin{equation}\label{ex1dof}
\begin{aligned}
&x_1'=x_2,\\
&x_2'=-0.0001[x_2^2+(x_1-1)^2-(1+\mu)^2]x_2-x_1+1,\\
&\Delta x_2|_{x\in \Gamma}= -( 1+Rx_2+\mu x_2^3)x_2+\mu^2,
\end{aligned}
\end{equation}
where  $R=0.9$  and  $\Gamma=\{x| x_1=0, x_2\leq0 \}.$ It is easy to see that system (\ref{ex1dof}) is of the form (\ref{smallpar}) and $\Phi(x_1,x_2)=x_1=0.$ The system  has a  periodic solution \begin{equation}\label{periodic1}\Psi_\mu(t)=(1+(1+\mu)\cos(t),-(1+\mu)\sin(t)),\end{equation} where $t\in\mathbb{R}$ for $\mu\in(-2,0].$ 

The   generating   system of (\ref{ex1dof})   has   the following   form

\begin{equation}\label{ex1dof1}
\begin{aligned}
&x_1'=x_2,\\
&x_2'=-0.0001[x_2^2+(x_1-1)^2-1]x_2-x_1+1,\\
&\Delta x_2|_{x\in \Gamma}= -(1+Rx_2)x_2,
\end{aligned}
\end{equation} 
  and  admits   the  periodic solution $\Psi_0(t)= (1+\cos(t),-\sin(t)).$  By means of the equality $\langle \nabla\Phi(x^*),f(x^*)\rangle=\langle (1,0),(0,1)\rangle=0$ with  $x^*=(0,0)\in\partial\Gamma,$  it is easy to say  that $x^*$ is a grazing point of  $\Psi_0(t).$

Let us start with   the linearization of system \eqref{ex1dof1} around the periodic solution $\Psi_0(t).$  Consider a near solution $y(t)=y(t,0,y^*+\Delta y),$ where $\Delta y=(\Delta y_1,\Delta y_2),$ to the periodic solution $\Psi_0(t).$   Assume that $y(t)$  satisfies condition $(N1),$  and  it meets the surface of discontinuity $\Gamma$ at the moment $t=\xi$ and  at the point $\bar y=y(\xi,0,y^*+\Delta y).$ 
Considering the formula \eqref{dert} for the transversal point $\bar y=(\bar y_1,\bar y_2),$ the first component $\displaystyle{\frac{\partial \tau(\bar y)}{\partial y_1^0}}$ can be evaluated as $\displaystyle{\frac{\partial \tau(\bar y)}{\partial y_1^0}=-\frac{1}{\bar y_2}}.$ From the last equality, the singularity is seen at the grazing point.   By taking into account \eqref{derW_1} with \eqref{ex1dof1} and  $\displaystyle{\frac{\partial \tau(\bar y)}{\partial y_1^0}},$  we obtain that

\begin{eqnarray}\label{derW_44}
\displaystyle{\frac{\partial W_i(\bar y)}{\partial y_1^0}}= \begin{bmatrix} R\bar y_2-1 \\ -0.0001R(\bar y_2^2+(\bar y_1-1)^2-1)-2R(0.0001(\bar y_2^2+(\bar y_1-1)^2-1))\end{bmatrix}.
\end{eqnarray}

Similarly, taking into account  the formula \eqref{derW_5},   one   can evaluate   that  $\displaystyle{\frac{\partial \tau(\bar y)}{\partial y_2^0}=0.}$ This and  formula \eqref{derW_6}  imply  

\begin{eqnarray}\label{derW_77}
\displaystyle{\frac{\partial W_i(\bar y)}{\partial y_2^0}}= \begin{bmatrix} 0 \\ -2R\bar y_2\end{bmatrix}.
\end{eqnarray}

Joining \eqref{derW_44} and \eqref{derW_77}, the matrix $W_{iy}(\bar y)$ can be obtained as 
\begin{eqnarray}\label{derW_88}
\displaystyle{ W_{iy}(\bar y)}= \begin{bmatrix} R\bar y_2-1  & 0 \\ -0.0001R(\bar y_2^2+(\bar y_1-1)^2-1)-2R(0.0001(\bar y_2^2+(\bar y_1-1)^2-1))&-2R\bar y_2\end{bmatrix}.
\end{eqnarray}

The    last  expression  implies   continuity   of the  partial   derivatives   near the grazing point.   This validates condition   $(A2).$    

Then, evaluating  the  matrix in \eqref{derW_88} at $\bar y= y^*=(0,0),$ it is easy to obtain  

\begin{eqnarray}\label{derW_99}
\displaystyle{ W_{iy}(y^*)}= \begin{bmatrix}  -1  & 0 \\ 0.0003R &0\end{bmatrix},
\end{eqnarray}

and 

\begin{eqnarray} \label{betaexi} B_i=\begin{cases}  O_{2}, \quad &\mbox{if } \quad \mbox{ $(N1)$ is valid,} \\ 
 \begin{bmatrix}  -1  & 0 \\ 0.0003R &0\end{bmatrix}, \quad &\mbox{if } \quad \mbox{ $(N2)$ is valid.}  \end{cases}
\end{eqnarray}

To  verify condition $(A3),$    let   us specify     the region 
$$H=\{(y_1,y_2)|y_2<\sqrt{1-(y_1-1)^2},\, 0\leq y_1\leq 1\}.$$ 
For the grazing solution $\Psi_0(t),$ we have that $\theta_{i+1}-\theta_i=2\pi.$ Consider a near solution $y(t)=(y_1(t),y_2(t))=y(t,0,\bar y)$ to $\Psi(t).$ To satisfy the condition $(N1),$ take $\bar y=(\bar y_1,\bar y_2)\in H.$ The orbit of $y(t)$ is below  the grazing orbit. Fix points $y=(y_1,y_2) \in H$ and $\psi=(\psi_1,\psi_2)$ of the orbits $y(t)$ and $\Psi_0(t),$ respectively such that $0\leq y_1=\psi_1\leq 1$ and $\psi_2<0.$ Since of the equation $y_1^\prime=y_2,$ the speed of $y_1(t)$ at $(y_1,y_2)$ is larger than the speed of $\Psi_1(t) $ at $(\psi_1,\psi_2).$ Consequently, one can find that $\tau(y)\leq \frac{\pi}{4}<2\pi$ for $y\in H.$ Thus, the condition $(A3)$ is valid and  Lemma \ref{conttau}  verifies the condition  $(C).$

 It is easy to demonstrate that the condition  $(A4)$ is valid such that  near solutions  to   the grazing one  are either continuous or discontinuous. That is, they don't intersect the line of discontinuity $\Gamma$ or intersect it permanently near to the   grazing point  and by means of the formula \eqref{betaexi}, the linearization system for \eqref{ex1dof1} around the grazing cycle $\Psi_0(t)$ consists of the following two  subsystems  
 \begin{equation}\label{ex1dofvar21}
\begin{aligned}
&u_1'=u_2,\\
&u_2'=-0.0001\sin(2t)u_1+0.0002 \sin^2(t)u_2,
\end{aligned}
\end{equation}
and 
 \begin{equation}\label{ex1dofvar22}
\begin{aligned}
&u_1'=u_2,\\
&u_2'=-0.0001\sin(2t)u_1+0.0002 \sin^2(t)u_2,\\
 &\Delta u|_{2\pi i}=\begin{bmatrix}  -1  & 0 \\ 0.0003R &0\end{bmatrix}u.
\end{aligned}
\end{equation}
The system \eqref{ex1dofvar21} + \eqref{ex1dofvar22}  is $(2\pi,1)$ periodic.    The Floquet multipliers of system \eqref{ex1dofvar21} + \eqref{ex1dofvar22}   are $\rho_1^{(1)}=1,$  $\rho_2^{(1)}=0.939,$ $\rho_1^{(2)}=1,$  $\rho_2^{(2)}=0.912.$ Thus, condition $(A5)$ is validated. Moreover, the conditions $(C1)-(C7)$ and $(A1),(A2)$ can be verified utilizing similar way  presented in  Example \ref{exdds}. Consequently,  Theorem \ref{thmorbitalstability} authenticates that the grazing periodic solution (cycle), $\Psi_0(t)$ of the system \eqref{ex1dof1} is orbitally  stable. The simulation results  demonstrating the orbital stability of $\Psi_0(t)$  are depicted in Figure \ref{ivanov}.

\begin{figure}[ht] 
\centering
\includegraphics[width=9 cm]{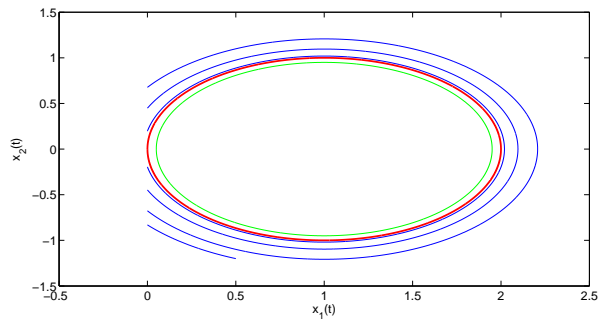}
\caption{The grazing   cycle  of system (\ref{ex1dof1})  is in red.   The blue arcs are the trajectory of the system  with initial   point $(0.5, 1.2)$ and the green   continuous   orbit  is  with initial value $(0.1,0).$  They   demonstrate   stability   of the grazing   solution.}
\label{ivanov}
\end{figure}

Next, we will investigate two sorts of  periodic solutions of system (\ref{ex1dof}) with a period $\mathscr{T}$ near to $2\pi.$   The first one is continuous and the second admits discontinuities once on a period.  For those solutions,  corresponding  linearization systems around the grazing cycle $\Psi_0(t)$ are  \eqref{ex1dofvar21} and  \eqref{ex1dofvar22}, respectively. Let us start with the continuous periodic solutions of (\ref{ex1dof}). For continuous periodic solution, we will consider the linearization system \eqref{ex1dofvar21}.

 To apply Theorem \ref{theosmall1}, denote $\Psi_0(0)= (\zeta_1^0,0). $ That   is, consider  $\zeta_2^0=0.$ Then, applying the above discussion, obtain that the Poincar$\grave{e}$ condition admits the form   of the  following   equations, 

\begin{equation}\label{example1}
\begin{aligned}
&\mathscr{S}_1(\mathscr{T},\zeta_1,\mu)=x_1(\mathscr{T},\zeta_1,\mu)-x_1=0,\\
&\mathscr{S}_2(\mathscr{T},\zeta_1,\mu)=x_2(\mathscr{T},\zeta_1,\mu) =0.
\end{aligned}
\end{equation} 

Because  solutions of  the system  \eqref{ex1dof1} have continuous  derivatives with respect to the time, phase variables and parameters,  we   can calculate the following determinant 
\begin{equation}\label{detex1}
\begin{vmatrix} \displaystyle{\frac{\partial \mathscr{S}_1(\omega,\zeta_1^0,0)}{\partial \mathscr{T}}} & \displaystyle{\frac{\partial \mathscr{S}_1(\omega,\zeta_1^0,0)}{\partial x_1^0} }\\ \displaystyle{\frac{\partial \mathscr{S}_2(\omega,\zeta_1^0,0)}{\partial \mathscr{T}}} & \displaystyle{\frac{\partial \mathscr{S}_2(\omega,\zeta_1^0,0)}{\partial x_1^0}}
\end{vmatrix}.
\end{equation}

First, we need the   monodromy matrix of the system \eqref{ex1dofvar21}. It is
 \begin{equation}\label{monodromy1}
\begin{bmatrix} 0.939 & -0.0001407 \\ -0.0003165 &1
\end{bmatrix}.
\end{equation}

It is easy to see that first column of the determinant (\ref{detex1})  is computed by utilizing (\ref{ex1dof1}) and the second column is evaluated by means of the first column of the matrix (\ref{monodromy1}). From this discussion, one can obtain that the   determinant  \eqref{detex1}  is equal  to

\begin{equation}\label{last1}
\begin{vmatrix} 0 & -0.061\\
                1 &  -0.0003165
\end{vmatrix}=0.061\neq 0.
\end{equation}

  Thus, in the light of Theorem \ref{theosmall1},  we can conclude that for sufficiently small  $|\mu|$ there exists a  unique periodic solution  of the system 
\begin{equation}\label{ex1dofode}
\begin{aligned}
&x_1'=x_2,\\
&x_2'=-0.0001[x_2^2+(x_1-1)^2-(1+\mu)^2]x_2-x_1+1.
\end{aligned}
\end{equation} 
It is exactly the cycle (\ref{periodic1}) with a period $\mathscr{T}=2\pi.$   If $\mu < 0,$   the    solution is separated from the   set  $\Gamma.$     Consequently, it  is a periodic   continuous solution of the equation (\ref{ex1dof}).    It  is   orbitally   stable by   the   theorem   for   continuous dynamics  \cite{Hirsh-smale},  since   of the   continuous   dependence of multipliers   on the   parameter. 
The   function $\Psi_{\mu}(t), \mu  >0,$  intersects $\Gamma$  and can not   be a solution  of    equation (\ref{ex1dof}).  Thus, the system  does not   admit  a continuous periodic solution near  to  
$\Psi_{0}(t),$  if the   parameter is positive.

Considering those solutions which have one moment of discontinuity in a period, one can find that the corresponding linearization of $\Psi_0(t)$ is the system \eqref{ex1dofvar22}.

The monodromy matrix of \eqref{ex1dofvar22} can be evaluated as   
  \begin{equation}\label{monodromy2}
\begin{bmatrix} 0.939 & -0.00052 \\ -0.000427 &1
\end{bmatrix}.
\end{equation} 
  It can be easily observed that the discontinuous solution intersects the line $x_1=0.$ For this reason, one can specify the first coordinate of the initial value as $\zeta_1=\zeta_1^0\equiv0.$  In the light of these discussions and the formula \eqref{periodpar},  the following equations are obtained:
\begin{equation}\label{example1small11}
\begin{aligned}
& \mathscr{S}_1(\mathscr{T},0,\zeta_2,\mu)=x_1(\mathscr{T},0,\zeta_2,\mu)=0, \\
&\mathscr{S}_2(\mathscr{T},0,\zeta_2,\mu)=x_2(\mathscr{T},0,\zeta_2,\mu)-\zeta_2=0.
\end{aligned}
\end{equation} 
Then, taking the derivative of the system \eqref{example1small11} with respect to $\mathscr{T}$ and $\zeta_2,$ and calculating it at $\mathscr{T}=\omega,$ $\zeta_2=\zeta_2^0=0,$ and for $\mu=0,$ the following determinant is obtained 
\begin{eqnarray}\label{smallexoneoneone111}
\begin{vmatrix} \displaystyle{\frac{\partial \mathscr{S}_1(\omega,0,\zeta_2^0,0)}{\partial \mathscr{T} }}&  \displaystyle{\frac{\partial \mathscr{S}_1(\omega,0,\zeta_2^0,0)}{\partial \zeta_2 }}\\ \displaystyle{\frac{\partial \mathscr{S}_2(\omega,0,\zeta_2^0,0)}{\partial \mathscr{T} }}&  \displaystyle{\frac{\partial \mathscr{S}_2(\omega,0,\zeta_2^0,0)}{\partial \zeta_2 }} \end{vmatrix} =
\begin{vmatrix}  0 & -0.0006\\ 1 & 0.0009\end{vmatrix}=-0.0006\neq 0.
\end{eqnarray}
Thus, condition $(A6)$ holds. Then,  utilizing Theorem \ref{theosmall1},  it is easy to conclude that for sufficiently small $\mu$ there exists a  unique periodic solution of the system (\ref{ex1dof}) with a period $\approx 2\pi.$  It is true that for positive as well as negative $\mu.$   Moreover, these solutions are orbitally asymptotically stable because of the continuous dependence of solutions on parameter and initial values  and they   meet the   discontinuity   line   transversally.

For  each   fixed  $\mu\neq 0,$  solutions near   to  the periodic   ones  intersect the line of discontinuity $\Gamma$ transversally  once   during the   time approximately  equal  to  the   period.   That is,   the smoothness which  is requested  for   the   application of the Poincar$\grave{e}$  condition   is  valid, since   the   smoothness   for  the grazing point  has already  been    verified. 
It  is   clear   that   there   can  not  be   another   solutions   with   period  close   to  $ 2\pi.$ 
Thus, one   can make the   following conclusion. The   original    system   (\ref{ex1dof})  admits   two    orbitally    stable   periodic  solutions,   continuous and discontinuous, if $\mu  < 0.$
There   is a   single   orbitally  stable continuous   solution  (grazing)   if $\mu = 0.$ Additionally, there  is a unique  discontinuous orbitally  stable   periodic solution for  positive values of  the parameter.
Consequently,   grazing  bifurcation of   cycles   appears  for  the system with  small parameter.

We have obtained  regular behavior in dynamics near grazing orbits by  the Poincar$\grave{e}$ small parameter analysis. Nevertheless,  outside the  attractors irregular phenomena  may be observed.

In Figure \ref{discont}, the solutions of the system (\ref{ex1dof}) with parameter $\mu=-0.2$ are depicted through simulations.  The red arcs are the trajectory of the system (\ref{ex1dof}) with initial value $(0.7, 0.05)$ and  the blue arcs are the trajectory of the system (\ref{ex1dof}) with initial value $(0.4, 0.05).$ It is seen that both red and blue trajectories approach the  discontinuous periodic solution of (\ref{ex1dof}), as time increases. So, the discontinuous cycle is orbitally stable trajectory. Moreover,  the green one is a continuous periodic trajectory of (\ref{ex1dof}) with initial value $(0, 0.05)$ and it is orbitally asymptotically stable. To sum up, there exists two  periodic solutions  of (\ref{ex1dof})  for the parameter $\mu=-0.2,$ one is continuous, the other one is discontinuous and  both solutions are orbitally asymptotically stable.  

\begin{figure}[ht] 
\centering
\includegraphics[width=9 cm]{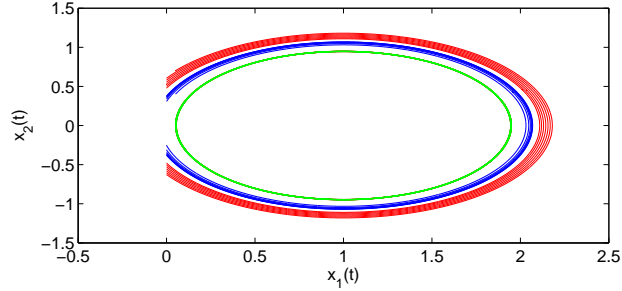}
\caption{The blue, red and green arcs constitute the trajectories of system (\ref{ex1dof}) with  $\mu=-0.2.$ The   first  two    approach   as time increases to  the discontinuous limit  cycle and the   third  one  is  the continuous limit  cycle itself. }
\label{discont}
\end{figure}
In Fig.  \ref{discontpert}, the red arcs are the orbit  of the system   with initial value $(0, 0.1)$    and the blue arcs are the trajectory of it with initial value $(0, 0.4).$ Both trajectories  approach to the discontinuous cycle of system   (\ref{ex1dof}),  as time increases. Thus,     Fig.  \ref{discontpert}   illustrates the existence of   the orbitally stable discontinuous periodic solution  if  $\mu=0.2.$  
\begin{figure}[ht] 
\centering
\includegraphics[width=10 cm]{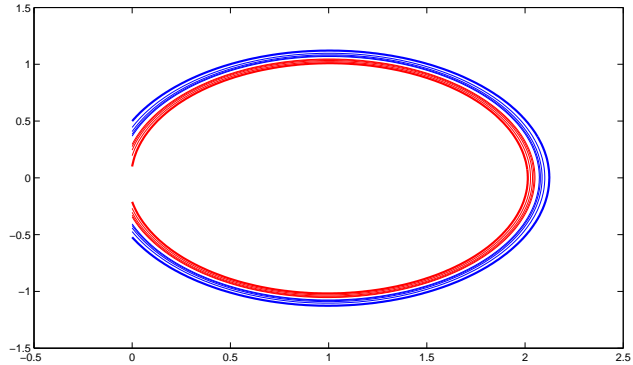}
\caption{The red and blue arcs constitute the trajectories  of the system (\ref{ex1dof}) with $\mu=0.2.$   Both orbits  approach to the discontinuous limit cycle, as time increases.}
\label{discontpert}
\end{figure}
\end{example}

\section{Conclusion}
In literature, the dynamics in the neighborhood of the grazing points \cite{budd2001}-\cite{ref24}, \cite{nusse94}, \cite{feigin70}, \cite{feigin78}, \cite{Nord1}-\cite{Nord2006} is generally analyzed through maps of the Poincar$\acute{e}$  type. The main analysis is conducted on complex dynamics behavior such as chaos and bifurcation \cite{budd2001}-\cite{ref24}, \cite{nusse94}, \cite{donde-hiskens}, \cite{feigin70}, \cite{Nord1}-\cite{Nord2006}. However,   there  is still no   sufficient  conditions for  the discontinuous motion to  admit   main features of dynamical  systems : the group property,  continuous and differentiable dependence on initial data and  continuation of motions, which are useful for both local and global analysis. Variational systems for grazing solutions have not been considered  in general  as well as  orbital stability theorem   and regular  perturbation theory    around cycles, despite, particular cases can be found in specialized papers.  See, for  example,  \cite{ref38}.  To investigate   these  problems  in the present  paper,  we  have applied   the  method of $B-$ equivalence  and  results on discontinuous dynamics developed  and    summarized   in \cite{ref1}.  In   our    analysis   the  grazing   singularity   is observed through  the gradient   of   the time function $\tau(x),$  since some of its coordinates are  infinite.  We  have found the components of the discontinuous  dynamical system  that is the vector field, surfaces of discontinuity and the equations of jump such that interacting they neutralize the effect of singularity. Then,  we  linearize   the system  at   the   grazing   moments and   this   brings the dynamics to regular analysis and   make suitable for the application. By means of the linearization, the theory can be understood as a part of the general theory of discontinuous dynamical system.  Thus, we have considered grazing phenomena as a   subject   of the general theory of   discontinuous dynamical  systems \cite{ref1},  discovered  a partition of set of solutions   near  grazing  solution such that we  determine linearization around a grazing solution is a collection of several linear impulsive systems with fixed moments of impulses. This constitutes the main novelty of the present paper. To linearize a solution around the grazing one, a system from the collection  is to be utilized.  This   result  has been  applied  to  prove  the orbital  stability   theorem.    The way of analysis in  \cite{ref1}-\cite{MAkhmethod}  continues in the present paper and it admits all attributes   which    are proper for continuous dynamics \cite{Hirsh-smale}. That is why,   we    believe that   the method can  be   extended   for   introduction   and  research   of graziness  in other types of dynamics  such as partial  and   functional differential   equations and others.   Next, we plan to apply the present  results and the method of investigation for problems initiated in \cite{budd2001}-\cite{ref24}, \cite{Nord97}-\cite{Nord2006}, \cite{Piir-virgin-Champhneys}.

\bibliographystyle{model1-num-names}
\bibliography{elsarticle-num.bst}

\end{document}